\documentclass[12pt]{article}
\usepackage{amsmath,amsfonts,amssymb,amsthm,amscd}

\catcode`\@=11
\renewcommand\section{\@startsection{section}{1}{\z@}%
                                      {-3.25ex\@plus -1ex \@minus -.2ex}%
                                      {1.5ex \@plus .2ex}%
                                      {\normalfont\large\bfseries}}
\catcode`\@=\active
\catcode`\@=11
\renewcommand\subsection{\@startsection{subsection}{2}{\z@}%
                                      {-3.25ex\@plus -1ex \@minus -.2ex}%
                                      {1.5ex \@plus .2ex}%
                                      {\normalfont\large\bfseries}}
\catcode`\@=\active

\newcommand{\w}{\omega}
\newcommand{\1}{{\bf 1}}
\newcommand{\id}{{\rm id}}
\newcommand{\Hom}{{\rm Hom\,}}

\DeclareMathOperator{\Aut}{Aut\,}
\DeclareMathOperator{\gr}{gr}

\newcommand{\End}{{\rm End\,}}

\DeclareMathOperator{\sgn}{sgn}
\DeclareMathOperator{\soc}{soc}
\DeclareMathOperator{\tr}{tr}

\newcommand\haru[2]{{\rm span}\{\,#1\,|\,#2\,\}}

\newcommand\Z{\mathbb{Z}}
\newcommand\Zpos{\Z_{\geq0}}
\newcommand\Zplus{\Z_{>0}}

\newcommand\C{\mathbb{C}}

\newcommand\h{\mathfrak{h}}
\newcommand\g{\mathfrak{g}}

\newcommand{\NO}{\,{\raise0.25em\hbox{$\mathop{\hphantom {\cdot}}\limits^{_{\circ}}_{^{\circ}}$}}\,}

\newcommand\spl{\mathfrak{sl}_2(\C)}
\newcommand\filt{\mathcal{L}}
\newcommand\trip{{SF}}
\newcommand\tripeven{\trip^+}
\newcommand\tripodd{\trip^-}
\newcommand\triptw{\trip(\theta)}
\newcommand\triptweo{\trip(\theta)^{\pm}}
\newcommand\triptweven{\trip(\theta)^{+}}
\newcommand\triptwodd{\trip(\theta)^{-}}
\newcommand\T{\mathcal{T}}
\newcommand\Teven{\T^{+}}

\newcommand\Ttw{\T(\theta)}

\newcommand\indetrip{\widehat{\trip}}

\newtheorem{theorem}{Theorem}[section]
\newtheorem{proposition}[theorem]{Proposition}
\newtheorem{lemma}[theorem]{Lemma}
\newtheorem{corollary}[theorem]{Corollary}
\theoremstyle{definition}

\theoremstyle{remark}
\newtheorem{remark}[theorem]{\bf Remark}

\numberwithin{equation}{section}

\makeatletter
\@addtoreset{equation}{section}

\makeatother

\begin{document}
\begin{center}
{\large\bf A $\Z_2$-orbifold model of \\
the symplectic fermionic vertex operator superalgebra}
\end{center}
\begin{center}
{\bf Toshiyuki Abe\footnote{This research is supported in part by a grant from Japan Society for the Promotion of Science}}\\
abe@math.sci.ehime-u.ac.jp
\vskip1ex
Mathematical Sciences, Faculty of Science, Ehime university\\
2-5 Bunkyocho, Matsuyama, Ehime 970-8577
\end{center}
\begin{abstract}
We give an example of an irrational $C_2$-cofinite vertex operator algebra whose central charge is $-2d$ for any positive integer $d$. 
This vertex operator algebra is given as the even part of the vertex operator superalgebra generated by $d$ pairs of symplectic fermions, and it is just the realization of the $c=-2$-triplet algebra given by Kausch in the case $d=1$.  
We also classify irreducible modules for this vertex operator algebra and determine its automorphism group.   
\end{abstract}
\tableofcontents

\section{Introduction}
Vertex operator algebras have been studied for twenty years. 
One of the motivations for their introduction comes from the conformal field theory, and many mathematicians and physicists are studying conformal field theory from the view point of vertex operator algebras. 
It is expected that the theory of vertex operator algebras satisfying a finiteness condition called $C_2$-cofiniteness corresponds to the theory of rational conformal field theories and that an irrational $C_2$-cofinite vertex operator algebra, that is a $C_2$-cofinite one which admits reducible indecomposable modules, gives a model of a logarithmic rational conformal field theory.     
In this paper, we exhibit an example of an irrational $C_2$-cofinite vertex operator algebra of central charge $-2d$ for any positive integer $d$.
We also classify its irreducible modules and determine its automorphism group.    

The notion of $C_2$-cofiniteness is quite important in the representation theory of vertex operator algebras.   
It was introduced by Y. Zhu in \cite{Zh1} as a sufficient condition for the existence of differential equations which one point functions on the torus should satisfy. 
This technical condition becomes essential to study the representation theory of vertex operator algebras because it enables us to prove many properties which rational conformal field theories have. 
For example, a $C_2$-cofinite vertex operator algebra has a rational central charge and rational conformal weights, and the number of its irreducible modules and fusion rules among them are finite (cf. \cite{Miya} and references in there). 
It is believed that the theory of $C_2$-cofinite vertex operator algebras corresponds to rational conformal field theory. 
The rationality of a vertex operator algebra, which is the property that any module equipped with a lower truncated $\Z$-grading is completely reducible, plays a special role when the vertex operator algebra is $C_2$-cofinite. 
In this case, the factorization property of the space of conformal blocks on the projective line (see \cite{NT}) and the Verlinde formula (see \cite{Huang1}) hold. 

The triplet algebras are examples of models of logarithmic rational conformal field theories. 
They are current algebras given as extensions of the Virasoro algebra of central charges $c=c_{p,1}=1-6(p-1)^2/p$ for integers $p\geq2$ and generated by the stress energy tensor and $3$ primary fields of weight $2p-1$. 
In the case $p=2$, the representation theory of the triplet algebra is studied in \cite{K1} and \cite{GaKa}, where the existence of reducible indecomposable modules is shown. 
A realization of the $c=-2$ triplet algebra by means of a pair of symplectic fermions is given in \cite{K1}. 
This realization gives a hint to construct a vertex operator algebra with the central charge $-2d$ for any positive integer $d$ such that it corresponds to the $c=-2$ triplet algebra in the case $d=1$. 

We briefly explain the construction of our vertex operator algebras. 
It is quite similar to that of the $\Z_2$-orbifold model of the free bosonic vertex operator algebra (see \cite{DN1} and \cite{DN3}). 
First we consider a finite dimensional vector space $\h$ with a nondegenerate \textit{skew-symmetric} bilinear form $\langle\,\cdot\,,\cdot\,\rangle$. 
Then $d=\dim\h/2$ is a nonnegative integer. 
We next consider the nontrivial central extension of the supercommutative Lie superalgebra $\h\otimes\C[t^{\pm1}]$ by the one dimensional center $\C K$ and a $2$-cocycle associated to $\langle\,\cdot\,,\cdot\,\rangle$.
The Fock space $\trip$ generated from the highest weight vector $\1$ characterized by the properties that $\h\otimes\C[t].\1=0$ and $K.\1=\1$ naturally has a structure of a vertex superalgebra such that the vector $\1$ is its vacuum vector.
In this vertex superalgebra, we consider the Virasoro vector given by $\sum_{i=1}^{d}(e^i\otimes t^{-1})(f^i\otimes t^{-1})\1$, where $\{e^i,f^i\}_{1\leq i\leq d}$ is a basis of $\h$ such that $\langle e^i,e^j\rangle=\langle f^i,f^j\rangle=0$ and $\langle f^i,e^j\rangle=\delta_{i,j}$ for any $1\leq i,j\leq d$.   
Then $\trip$ becomes a vertex operator superalgebra of central charge $c=-2d$. 
This vertex operator superalgebra has the canonical involution $\theta$ associated to the $\Z_2$-grading. 
The fixed point set $\tripeven$ of $\trip$ by $\theta$, or the even part of $\trip$ is our vertex operator algebra. 

The odd part of $\trip$ is an irreducible $\tripeven$-module. 
Other two irreducible modules appear in the $\theta$-twisted $\trip$-module. 
We show that any irreducible $\tripeven$-module is isomorphic to one of the irreducible modules above by investigating the structure of Zhu's algebra of $\tripeven$ introduced in \cite{Zh1}. 
By construction, the vertex operator algebra $\tripeven$ with the central charge $-2d$ contains a tensor product of $d$ copies of $\tripeven$ with $d=1$ as a full vertex operator subalgebra. 
Hence the $C_2$-cofiniteness of $\tripeven$ follows from that of $\tripeven$ with $d=1$. 
The $C_2$-cofiniteness of $\tripeven$ with $d=1$ can be shown by using the explicit forms of the relevant null vectors described in \cite{GaKa} via the realization (recently the $C_2$-cofiniteness of the $c_{p,1}$-triplet algebra with $p\geq2$ is proved in \cite{CF}). 
But in this article, we prove the $C_2$-cofiniteness of $\tripeven$ with $d=1$ directly. 

We also construct two reducible indecomposable $\tripeven$-modules. 
They can be obtained as $\tripeven$-submodules of the $\trip$-module generated from a singular vector $\hat{\1}$ characterized by the properties that $\h\otimes t\C[t].\hat{\1}=0$ and $K.\hat{\1}=\hat{\1}$. 
The existence of such modules shows that $\tripeven$ is not rational.       
From the construction of irreducible modules, we can express the all irreducible characters for $\tripeven$ by means of the Dedekind eta function. 
As is known for the case of the triplet algebra, we find that there appear polynomials of the logarithmic terms $2\pi i\tau=\log q$ of the coefficients of irreducible characters in the modular transformation for the transformation $\tau\rightarrow-\tau^{-1}$. 
The modular invariance for such a vertex operator algebra is studied in \cite{Miya}.
We expect that suitable subquotients of the reducible indecomposable $\tripeven$-modules above are interlocked by some symmetric linear function of the first Zhu's algebra $A_1(\tripeven)$, and that the space of $1$-point functions on the torus may be spanned by the pseudotrace functions associated with these modules and irreducible modules included in the $\theta$-twisted module.     
We also determine the automorphism group of $\tripeven$ and show that it is isomorphic to $Sp(2d,\C)/\langle\pm1\rangle$.
Any element of $\Aut(\tripeven)$ is given by a natural lifting of a linear isomorphism of $\h$ which preserves the skew-symmetric bilinear form $\langle\,\cdot\,,\cdot\,\rangle$. 

This article is organized as follows.
In Section \ref{Secta}, we recall some definitions and results in the representation theory of vertex operator (super)algebras. 
In Section \ref{Secta1}, we review shortly a definition of vertex operator superalgebras, their modules and their automorphisms.
The definition of $C_2$-cofiniteness and Zhu's algebras are stated in Section \ref{Secta2}, where we give a relation between a set of generators of a vertex operator algebra and that of its Zhu's algebra.  
In Section \ref{Sectb}, we construct the vertex operator algebra $\tripeven$ and study its structure.
In Section \ref{Sectb1}, we give a construction of $\tripeven$, and give a suitable generating set of $\tripeven$ in Section \ref{Sectb2}. 
A proof of $C_2$-cofiniteness of $\tripeven$ is given in Section \ref{Sectb3}.
Section \ref{Sectc} is devoted to the classification of irreducible $\tripeven$-modules. 
The way is similar to that of the vertex operator algebras $M(1)^+$ and $V_L^+$ for positive definite even lattices $L$ (see \cite{DN1}--\cite{DN3} and \cite{AD}). 
In Section \ref{Sectc1} we give a construction of an irreducible $\theta$-twisted $\trip$-module.
In Section \ref{Main}, the main theorem is stated. 
Section \ref{Sectc2} and Section \ref{Sectc3} are proofs of the main theorem for $d=1$ and $d>1$ respectively. 
In Section \ref{Sectc4} we construct reducible indecomposable $\tripeven$-modules which proves that $\tripeven$ is irrational. 
We calculate irreducible characters and determine the automorphism group of $\tripeven$ in Section \ref{Secte2}. 

The author thanks K. Nagatomo and A. Matsuo for fruitful comments and insightful advice.    
He also thanks C. Dong, Y. Arike, for reading the manuscript and giving him some errors.       

\section{Preliminaries}\label{Secta}
In this section we recall some notions and results in the representation theory of vertex operator (super)algebras.  
Throughout the paper we use the notations $\Zpos$ and $\Zplus$ for the set of all nonnegative integers and positive integers respectively.  
\subsection{Vertex operator superalgebras, their modules and rationality}\label{Secta1}
A {\em vertex operator superalgebra} is a $4$-tuple $(V,Y(\,\cdot\,,z),\1,\w)$ which consists of a $\Z_2$-graded vector space $V=V^{\bar 0}\oplus V^{\bar 1}$, a linear map $Y(\,\cdot\,,z):V\rightarrow \End V[[z,z^{-1}]]$ which maps $a\in V$ to $Y(a,z)=\sum_{n\in\Z}a_{(n)}z^{-n-1}$ with $a_{(n)}\in\End V$ and vectors $\1,\w\in V^{\bar 0}$ satisfying the following axioms (1)--(6) (see \cite{L3}, \cite{Xu} or \cite{Kac}): 

(1)  
For $a\in V^{\bar{k}},b\in V^{\bar{l}}\,(k,l\in\{0,1\})$ and $n\in\Z$, $a_{(n)}b\in V^{\bar{k}+\bar{l}}$.
Furthermore, $a_{(n)}b=0$ if $n$ is sufficiently large.

(2) 
For any $a\in V^{\bar{k}},b\in V^{\bar{l}}\,(k,l\in\{0,1\})$ and $p,q,r\in\Z$, the identity called the Borcherds identity 
\begin{align}\label{superBorcherds}
\sum_{i=0}^{\infty}\binom{q}{i}(a_{(p+i)}b)_{(q+r-i)}
=\sum_{i=0}^{\infty}(-1)^i\binom{p}{i}(a_{(p+q-i)}b_{(r+i)}-(-1)^{p+kl}b_{(p+r-i)}a_{(q+i)})
\end{align}
holds in $\End V$.

(3)
The vector $\1$ called the {\em vacuum vector} satisfies that $\1_{(m)}=\delta_{m,-1} \id_V$, $a_{(-1)}\1=a$ and $a_{(n)}\1=0$ for any $a\in V$, $m\in\Z$ and $n\in\Zpos$.
 
(4)
The set of operators $\{L_n,\id_V\}_{n\in\Z}$ with $L_n:=\w_{(n+1)}$ gives a representation of the Virasoro algebra on $V$ of central charge $c_V\in \C$; that is, 
\[
[L_m,L_n]=(m-n)L_{m+n}+\frac{m^3-m}{12}\delta_{m+n,0}c_V\id_V
\] 
for any $m,n\in \Z$.
The vector $\w$ is called the {\em Virasoro vector}.

(5) 
For any $a\in V$, $L_{-1}a=a_{(-2)}\1$. 

(6) 
$V$ is decomposed into a direct sum of finite dimensional eigenspaces $V_n$ for $L_0$ of eigenvalues $n\in\frac{1}{2}\Zpos$ as
\[
V=\bigoplus_{n\in\frac{1}{2}\Zpos}V_n. 
\] 
The eigenvalues for $L_0$ are called {\em weights}. 
\vskip1ex 
We often refer $V$ to be a vertex operator superalgebra for simplicity.  
This definition of ``vertex operator superalgebra" is also known in the literature as ``conformal vertex superalgebra" (see for example \cite{Kac}).  
We see that axiom (5) is equivalent to the $L_{-1}$-derivative property $Y(L_{-1}a,z)=\frac{d}{dz}Y(a,z)$ for $a\in V$. 
Since $L_0\1=0$ and $L_0\w=2\w$, $\1\in V_0$ and $\w\in V_2$. 

We should note that $V^{\bar0}=\bigoplus_{i\in \frac{1}{2}\Z}V^{\bar0}\cap V_i$ and $V^{\bar1}=\bigoplus_{i\in \frac{1}{2}\Z}V^{\bar1}\cap V_i$ in general. 
In some literature, it is assumed that $V^{\bar{0}}=\bigoplus_{n\in\Zpos}V_n$ and $V^{\bar{1}}=\bigoplus_{n\in\frac{1}{2}+\Zpos}V_n$ in the definition of vertex operator algebras (cf. \cite{L3,Xu}). 
We do not assume here these conditions but the conditions  
\begin{align}\label{xxeo}
V_0=\C\1\quad\text{and } V_n=0\quad\text{for $n\in \frac{1}{2}+\Zpos$}. 
\end{align} 
  
A {\em vertex operator algebra} is a vertex operator superalgebra $V$ satisfying that $V^{\bar1}=0$ and $V_n=0$ for $n\in\frac{1}{2}+\Zpos$ (see \cite{FLM}, \cite{FHL}, \cite{MN} or \cite{LL}). 
By definition the following identity holds in $\End V$: 
\begin{align}\label{Borcherds}
\sum_{i=0}^{\infty}\binom{q}{i}(a_{(p+i)}b)_{(q+r-i)}
=\sum_{i=0}^{\infty}(-1)^i\binom{p}{i}(a_{(p+q-i)}b_{(r+i)}-(-1)^{p}b_{(p+r-i)}a_{(q+i)})
\end{align}
for any $a,b\in V$ and $p,q,r\in \Z$. 

An {\em automorphism} of a vertex operator superalgebra $V$ is a linear isomorphism $g$ of $V$ such that $g(\w)=\w$ and that $g(Y(a,z)b)=Y(g(a),z)g(b)$ for any $a,b\in V$. 
We denote by $\Aut V$ the group of all automorphisms of $V$. 
Axiom (1) shows that the map $\theta:V\rightarrow V$ defined by $\theta(a+b)=a-b$ for $a\in V^{\bar0}$ and $b\in V^{\bar 1}$ is an automorphism.

Let $V$ be a vertex operator superalgebra and $S$ a subset of $V$ consisting of eigenvectors for $L_0$ in $V^{\bar0}$ or $V^{\bar1}$. 
If $V$ is spanned by vectors of the form 
\[
u^1_{(-n_1)}u^{2}_{(-n_2)}\cdots u^{r}_{(-n_r)}\1
\]
for $u^i\in S$ and $n_i\in\Zplus$, then $V$ is called {\em strongly generated} by $S$ following \cite{Kac}. 
We call the subset $S$ a {\em set of field generators} of $V$. 
A vertex operator superalgebra $(U, Y_{U}(\,\cdot\,,z),\1_U,\w_U)$ is called a vertex operator subsuperalgebra of $V$ if $U\subset V$, $\1_U=\1$ and $Y_U(u,z)v=Y(u,z)v$ for any $u,v\in U$. 
If $\w_U=\w$ then $U$ is called a {\em full} vertex operator subsuperalgera $V$, or it is said that $U$ is {\em conformally embedded} in $V$. 

Let $V$ be a vertex operator superalgebra and $g$ an automorphism commuting with $\theta$ of order $T$.
Then $V$ is decomposed into a direct sum of the eigenspaces $V^{r}$ of $g$ with eigenvalue $e^{-\frac{2\pi i r}{T}}$ with $0\leq r\leq T-1$. 
A {\em weak $g$-twisted $V$-module} is a pair $(M,Y(\,\cdot\,,z))$ of a $\Z_2$-graded vector space $M=M^{\bar0}\oplus M^{\bar1}$ and a linear map $Y(\,\cdot\,,z):V\rightarrow \End M[[z^{\frac{1}{T}},z^{-\frac{1}{T}}]],a\mapsto Y(a,z)=\sum_{n\in\frac{1}{T}\Z}a_{(n)}z^{-n-1}$ satisfying the following axioms (1)--(4):

(1)  
For $a\in V^{\bar{k}},u\in M^{\bar{l}}\,(k,l\in\{0,1\})$ and $n\in\frac{1}{T}\Z$, $a_{(n)}u\in M^{\bar{k}+\bar{l}}$.
Furthermore,  $a_{(n)}u=0$ if $n$ is sufficiently large.

(2)  
If $a\in V^r$ for $0\leq r\leq T-1$ and $n\notin\frac{r}{T}+\Z$, then $a_{(n)}=0$ on $M$. 

(3) 
For any $a\in V^{\bar{k}}\cap V^s,b\in V^{\bar{l}}\cap V^t\,(k,l\in\{0,1\},0\leq s,t \leq T-1)$, $p\in\Z$ and $q\in\frac{s}{T}+\Z,r\in\frac{t}{T}+\Z$, the following identity holds in $\End M$:
\begin{align}\label{supertwBorcherds}
\sum_{i=0}^{\infty}\binom{q}{i}(a_{(p+i)}b)_{(q+r-i)}
=\sum_{i=0}^{\infty}(-1)^i\binom{p}{i}(a_{(p+q-i)}b_{(r+i)}-(-1)^{p+kl}b_{(p+r-i)}a_{(q+i)}).
\end{align}

(4)
The vacuum vector $\1$ satisfies $\1_{(n)}=\delta_{n,-1}\id_M$ for any $n\in\frac{1}{T}\Z$. 

\vskip1ex
If we set $L_n=\w_{(n+1)}$ for any $n\in\Z$, then $\{L_{n},\id_M\,|\,n\in\Z\}$ gives a representation of the Virasoro algebra on $M$ of central charge $c_V$. 
The $L_{-1}$-derivative property $Y(L_{-1}a,z)=\frac{d}{dz}Y(a,z)$ also holds for any $a\in V$ (see \cite{Xu} for example). 
  
Let $V$ be a vertex operator superalgebra satisfying \eqref{xxeo}. 
A weak $g$-twisted $V$-module is called \textit{$\frac{1}{T}\Zpos$-gradable} if $M$ has a $\frac{1}{T}\Zpos$-grading as $M=\bigoplus_{n\in\frac{1}{T}\Zpos}M{(n)}$ and satisfies the condition 
\[
a_{(n)}M{(m)}\subset M{(k+m-n-1)}
\]
for any $a\in V_k\,(k\in\Zpos)$ and $m,n\in\frac{1}{T}\Z$, where we set $M(n)=0$ for $n<0$. 
A vertex operator superalgebra is said to be {\em $g$-rational} if any $\frac{1}{T}\Zpos$-gradable weak $V$-module is completely reducible.

We call a weak $g$-twisted $V$-module $M$ a {\em $g$-twisted $V$-module} if $M$ is finite generated and if for any $u\in M$, the subspace spanned by vectors of the form $a_{(n)}u$ for $a\in V_k\,(k\in\Z_{\geq0})$ and $n\geq k-1$ is finite dimensional (cf.  \cite{NT}).
We see that any $V$-module $M$ is a direct sum of generalized eigenspaces for $L_0$ and that each generalized eigenspace is finite dimensional. 
We denote by $M_{(\lambda)}$ the generalized eigenspace for $L_0$ of eigenvalue $\lambda\in\C$ and say an eigenvalue for $L_0$ to be a weight. 
If $M=\bigoplus_{n=0}^\infty M_{(\lambda+n)}$ and $M_{(\lambda)}\neq 0$ for some $\lambda\in\C$ then we call the weight $\lambda$ the {\em lowest weight} of $M$. 
A $V$-module whose $V$-submodule is either $0$ or itself is called {\em irreducible}. 
If $V$ is irreducible as a $V$-module, $V$ is said to be {\em simple}.  

For a weak $g$-twisted $V$-module $M$, we set $\Omega(M)$ the subspace of $M$ which consists of the vectors $u$ such that $a_{(n)}u=0$ for any $a\in V_k\,(k\in\Zpos)$ and $n>k-1$. 
A vector in $\Omega(M)$ is called a {\em singular vector} in $M$. 
If $u\in M$ satisfies the condition that $L_1u=0$ then $u$ is called {\em quasi-primary} and if $u$ satisfies $L_nu=0$ for any $n\in\Zplus$ then called {\em primary}. 
We see that $u\in M$ is primary if and only if $L_1u=L_2u=0$.  
In the case $g=\id_V$, we refer simply a weak $g$-twisted module, a $g$-twisted $V$-module, $g$-rational to a weak module, a module and rational respectively. 
 
Let $V$ be a vertex operator superalgebra, $g$ an automorphism of finite order, $M$ a  $g$-twisted module and $N$ a $g^{-1}$-twisted $V$-module.  
Then a bilinear map $(\,\cdot\,,\cdot\,)$ from $M\times N$ to $\C$ is called {\em  invariant} if 
\begin{align}\label{kq;aaa}
(Y(a,z)u,v)=(-1)^{kl}(u,Y(e^{zL_{1}}(-z^{-2})^{L_0}a,z^{-1})v)
\end{align}
for any $a\in V^{\bar{k}}$, $u\in M^{\bar{l}}\,(k,l\in\{0,1\})$ and $v\in N$. 
\begin{proposition}\label{iwiaslka}
Let $M$ be a $g$-twisted $V$-module, $N$ a $g^{-1}$-twisted $V$-module and $(\,\cdot\,,\cdot\,)$ is a bilinear map from $M\times N$ to $\C$. 
Let $S$ be a set of field generators of $V$.
Suppose that $S$ consists of eigenvectors for $g$ and that for any $u\in M$, $(u,N_{(\lambda)})=0$ for any complex number $\lambda$ whose real part is sufficiently small.   
Then the bilinear map $(\,\cdot\,,\cdot\,)$ from $M\times N$ to $\C$ is invariant if and only if \eqref{kq;aaa} holds for any $a\in S$, $u\in M$ and $v\in N$. 
\end{proposition}
\begin{proof}
Set $D(N)$ to be the subspace of $N^*=\Hom_{\C}(N,\C)$ consisting of  vectors $f\in N^*$ satisfying that $f(Y(a,z)u)\in\C[z^{\pm\frac{1}{T}}]$ for any $a\in V$ and $v\in N$. 
Then $D(N)=D(N)^{\bar0}\oplus D(N)^{\bar1}$ becomes a weak $g$-twisted $V$-module with 
\[
(Y(a,z)f)(v)=(-1)^{kl}f(Y(e^{zL_{1}}(-z^{-2})^{L_0}a,z^{-1})v)
\]
for $a\in V^{\bar{k}}$, $f\in D(N)^{\bar{l}}$ and $v\in N$ (cf \cite{Li4}), where $D(N)^{\bar0}=D(N)\cap\Hom_{\C}(N^{\bar0},\C)$ and $D(N)^{\bar1}=D(N)\cap\Hom_{\C}(N^{\bar1},\C)$. 

We then have a linear map $\eta:M\rightarrow D(N),u\mapsto (u,\cdot\,)$ for $u\in M$.    
By the assumption, $\eta$ is well-defined and commutes with the action of $a_{(n)}$ for any $a\in S$ and $n\in\frac{1}{T}\Z$.
Since $V$ is strongly generated by $S$, $\eta$ commutes with the action of $a_{(n)}$ for any $a\in V$ and $n\in\frac{1}{T}\Z$.
Therefore, $\eta$ is a $V$-module homomorphism. 
This implies that \eqref{kq;aaa} holds for any $a\in V^{\bar{k}},u\in M^{\bar{l}}$ and $b\in N$.
The converse is clear. 
\end{proof}

\subsection{Zhu's algebra and $C_2$-cofiniteness condition}\label{Secta2}
In this section, we recall the notion of Zhu's algebra and $C_2$-cofiniteness for vertex operator algebras. 
 
Let $V$ be a vertex operator algebra. 
We consider the subspace $C_2(V)$ of $V$ spanned by vectors of the form $a_{(-2)}b$ for $a,b\in V$. 
Zhu found that the quotient space $V/C_2(V)$ has a commutative associative algebra structure in \cite{Zh1}.  
We denote by $\overline{a}$ the image of $a\in V$ in $V/C_2(V)$.
Then the product $\overline{a}\cdot\overline{b}$ of $\overline{a}$ and $\overline{b}$ for $a,b\in V$ is defined by 
\[
\overline{a}\cdot\overline{b}=\overline{a_{(-1)}b}.
\]
Since $a_{(-1)}b\equiv b_{(-1)}a\mod L_{-1}V$ for $a,b\in V$, $V/C_2(V)$ is commutative. 
We consider a set of field generators $S$ of $V$. 
Then one can see that $V/C_2(V)$ is spanned by vectors of the form 
\[
\overline{a^1}\cdot\overline{a^2}\cdot\cdots\cdot\overline{a^r}
\] 
with $a^i\in S$. 
Hence we have 
\begin{proposition}\label{iwiwsoi}
Let $S$ be a set of field generators of $V$. 
Then $V/C_2(V)$ is generated by $\bar{S}=\{\bar{a}\,|\,a\in S\}$ as a commutative algebra.  
\end{proposition}
 
Next we recall the definition of Zhu's algebra. 
Let $O(V)$ be the subspace of $V$ spanned by vectors of the form 
\[
\sum_{i=0}^{\infty}\binom{k}{i}a_{(i-2)}b
\]
for any vector $a\in V_k\,(k\in\Zpos)$ and $b\in V$, and set $A(V):=V/O(V)$. 
We write $[a]$ for the image of $a\in V$ in $A(V)$.  
We then have the well-known identity $[L_{-1}a]=-[L_0a]$ for any $a\in V$. 
 
In \cite{Zh1}, it is proved that $A(V)$ becomes an associative algebra with the product $[a]*[b]$ which is defined by 
\[
[a]*[b]=[a*b]=\sum_{i=0}^{\infty}\binom{k}{i}[a_{(i-1)}b]
\]
for $a\in V_k\,(k\in\Zpos)$ and $b\in V$. 
The linear isomorphism of $V$ which maps $a$ to $e^{L_1}(-1)^{L_0}a$ induces a linear isomorphism $\Phi$ of $A(V)$.
In fact, $\Phi$ is an anti-involution, i.e., $\Phi^2=\id$ and $\Phi([a]*[b])=\Phi([b])*\Phi([a])$ for any $a,b\in V$.
We note that if $a$ is a quasi-primary vector of weight $k$ then $\Phi([a])=(-1)^k[a]$. 
The involution $\Phi$ is useful to reduce some calculations when we find relations in $A(V)$. 
The following proposition is well known (see \cite{Zh1}). 
\begin{proposition}\label{loajda}
The image $[\1]$ is the unit of $A(V)$ and $[\w]$ is in the center of $A(V)$.
\end{proposition}
We often write $1$ for $[\1]$ in $A(V)$, and will use the fact that $[\1]$ is the unit and $[\w]$ in the center without referring this proposition.   

One of reasons to introduce the notion of Zhu's algebra is that the representation theory of Zhu's algebra  $A(V)$ is deeply related with that of the vertex operator algebra $V$.
\begin{theorem}\label{kiuwewytt} {\rm (\cite{Zh1})} 
Let $V$ be a vertex operator algebra. 
Then the following hold. 

(1) Let $M$ be a weak $V$-module.
Then the linear map $o:V\rightarrow \End \Omega(M),\,a\mapsto a_{(k-1)}$ for $a\in V_k\,(k\in\Z_{\geq0})$ induces a representation of $A(V)$ on $\Omega(M)$. 

(2) If $M$ is an irreducible $V$-module, then $\Omega(M)$ is irreducible as an $A(V)$-module. 

(3) For any irreducible $A(V)$-module $W$, there exists an irreducible $V$-module $M$ such that $\Omega(M)\cong W$ as $A(V)$-modules.

(4) The map $M\mapsto \Omega(M)$ induces a bijection from the set of inequivalent irreducible $V$-modules and that of inequivalent irreducible $A(V)$-modules.
\end{theorem} 

Let $S$ be a set of field generators of $V$. 
Following \cite[Proposition 3.3.2]{NT}, we show that $\{\,[a]\,|\,a\in S\,\}\subset A(V)$ is a generating set of $A(V)$ as an algebra. 
We consider a filtration $\{{F}_kA(V)\}_{k\in\Zpos}$ of $A(V)$ defined by 
\[
{F}_kA(V)=\left(\bigoplus_{i=0}^{k}V_i+O(V)\right)/O(V)
\] 
for any $k\in\Zpos$. 
Then we see that $F_kA(V)*F_lA(V)\subset F_{k+l}A(V)$ for any $k,l\in\Zpos$.
This implies that the multiplication of $A(V)$ induces an associative multiplication of the graded vector space $\gr_{\bullet} A(V)$, where 
\begin{align*}
\gr_{\bullet}A(V)=\bigoplus_{k=0}^\infty\gr_{k} A(V),\quad\gr_{k} A(V):={F}_kA(V)/{F}_{k-1}A(V)
\end{align*}
for any $k\in\Zpos$ and ${F}_{-1}A(V)=0$. 
Since $[a*b]-[b*a]\in F_{k+l-1}A(V)$ for $a\in V_k$ and $b\in V_l\,(k,l\in\Zpos)$, $\gr_{\bullet}A(V)$ is commutative. 

We consider a linear epimorphism $f:V\rightarrow \gr_{\bullet} A(V)$ defined by $f(a)=[a]+{F}_{k-1}A(V)\in \gr_{k} A(V)$ for $a\in V_k$. 
Since $f(a_{(-2)}b)=0$ for any $a,b\in V$, $f$ induces a linear epimorphism $\overline{f}$ from $V/C_2(V)$ to $\gr_{\bullet}A(V)$. 
In fact, we see that the epimorphism is of algebras because 
\[
\overline{f}(\overline{a}\cdot\overline{b})=f(a_{(-1)}b)=[a_{(-1)}b]+F^{k+l-1}A(V)= [a*b]+F^{k+l-1}A(V)=\overline{f}(\overline{a})\overline{f}(\overline{b})
\]
for $a\in V_k$ and $b\in V_l$. 
Now we have the following proposition. 
\begin{proposition}\label{eert}
Let $V$ be a vertex operator algebra and $S$ a set of field generators. 
Then $A(V)$ is generated by the set $\{[a]\,|\,a\in S\}$ as an associative algebra.  
\end{proposition}
\begin{proof}
Let $\bar{f}:V/C_2(V)\rightarrow \gr_{\bullet} A(V)$ be the algebra epimorphism above. 
By Proposition \ref{iwiwsoi}, we see that $\gr_{\bullet} A(V)$ is generated by $\bar{f}(\bar{a})=[a]+F^{k}A(V)$ for $k\in\Zpos$ and $a\in S\cap V_k$. 
Now by using induction on $k$, we can show that $F^{k}A(V)$ is contained in a subalgebra of $A(V)$ generated by all $[a]$ with $a\in S$ for any $k\in\Zpos$.    
\end{proof}

For a weak $V$-module $M$, we set $C_2(M)=\haru{a_{(-2)}u}{a\in V,u\in M}$. 
A weak $V$-module $M$ is called {\em $C_2$-cofinite} if $M/C_2(M)$ is finite dimensional.
We note that if $V$ is $C_2$-cofinite then Zhu's algebra $A(V)$ is of finite dimension. 
Therefore, there are only finitely many isomorphism classes of irreducible $V$-modules by Theorem \ref{kiuwewytt}. 
The following theorem is also one of the most remarkable results of the $C_2$-cofiniteness condition (see \cite{GN}, \cite{Buhl}, \cite{ABD}, \cite{NT} or \cite{Miya}). 
\begin{theorem}
Let $V$ be a $C_2$-cofinite vertex operator algebra. 
Then any finite generated weak $V$-module is a $C_2$-cofinite $V$-module. 
\end{theorem}
As a corollary, we have 
\begin{proposition}\label{rrlsiw}
Let $V$ be a vertex operator algebra and $U$ its full vertex operator subalgebra.
If $U$ is $C_2$-cofinite then $V$ is $C_2$-cofinite. 
\end{proposition} 

\section{The vertex operator algebra $\tripeven$}\label{Sectb}
In this section we construct the vertex operator superalgebra $\trip$ associated with a finite dimensional vector space $\h$. 
The even part $\tripeven$ of $\trip$ is the desired vertex operator algebra. 
We also find a set of field generators of $\tripeven$ and prove that $\tripeven$ is $C_2$-cofinite. 
\subsection{A construction of the vertex operator algebra $\tripeven$}\label{Sectb1}
Let $\h$ be a finite dimensional vector space with a skew-symmetric nondegenerate bilinear form $\langle\,\cdot\,,\cdot\,\rangle$.
Then the dimension of $\h$ is even and there is a basis $\{e^i,f^i\,|1\leq i\leq d\}$ such that 
\[
\langle e^i,e^j\rangle=\langle f^i,f^j\rangle=0\quad\text{and }\langle e^i,f^j\rangle=-\langle f^j,e^i\rangle=-\delta_{i,j}
\]
for any $1\leq i,j\leq d$, where we set 
\[
d=\frac{\dim \h}{2}. 
\] 
We call such a basis a {\em canonical basis} of $\h$ and denote it by $\{(e^i,f^i)\}_{1\leq i\leq d}$ or $\{(e^i,f^i)\}$ if the dimension of $\h$ is obvious. 

Now we consider the Heisenberg superalgebra $\hat{L}(\h):=\h\otimes\C[t,t^{-1}]\oplus\C K$ such that $\C K$ is the even part, $\h\otimes\C[t,t^{-1}]$ is the odd part and that the commutation relations are given by 
\[
[\psi\otimes t^m,\psi'\otimes t^n]_{+}=m\langle\psi,\psi'\rangle\delta_{m+n,0}K
\] 
for $\psi,\psi'\in\h,\,m,n\in\Z$ and $[K,\hat{L}(\h)]=0$. 
Let $\mathcal{A}$ be the quotient algebra of the universal enveloping algebra $U(\hat{L}(\h))$ by the two sided ideal generated by $K-1$. 
The $\Z_2$-grading of $\hat{L}(\h)$ naturally induces a $\Z_2$-grading on $\mathcal{A}$ as an algebra. 
We denote its even part and odd part by $\mathcal{A}^{\bar0}$ and $\mathcal{A}^{\bar1}$ respectively. 

We denote by $\psi{(m)}$ the operator of left multiplication by $\psi\otimes t^m$ on $\mathcal{A}$ for $\psi\in\h$ and $m\in\Z$. 
We set $\mathcal{A}_{\geq0}$ the left ideal generated by $\psi(m)1$ with $\psi\in\h,m\in\Zpos$ and consider the left $\mathcal{A}$-module $\trip=:\mathcal{A}/\mathcal{A}_{\geq0}$. 
We note that $\trip\cong\Lambda(\h\otimes t^{-1}\C[t^{-1}])$ as vector spaces. 
Since $\mathcal{A}_{\geq0}=\mathcal{A}_{\geq0}\cap\mathcal{A}^{\bar0}\oplus\mathcal{A}_{\geq0}\cap\mathcal{A}^{\bar1}$, we have $\trip=\trip^{\bar0}\oplus \trip^{\bar1}$, where $\trip^{\bar{i}}=\mathcal{A}^{\bar{i}}/(\mathcal{A}_{\geq0}\cap\mathcal{A}^{\bar{i}})$ for $\bar{i}\in \Z_2$. 

For any $n\in\Z$, we have a linear map 
\[\trip\rightarrow\End (\trip),a\mapsto a_{(n)}\]
defined by  
\begin{align}\label{lqqljq}
a_{(n)}=\sum_{\substack{i_j\in\Z\\\sum_{j=1}^ri_j=-\sum_{j=1}^rn_j+n+1}}\binom{-i_1-1}{n_1-1}\cdots\binom{-i_r-1}{n_r-1}\NO\psi^1(i_1)\cdots\psi^r(i_r)\NO
\end{align}
for $a=\psi^1(-n_1)\cdots\psi^r(-n_r)1$ with $\psi^i\in\h$ and $n_i\in\Zplus$, where $\binom{m}{n}=\frac{m(m-1)\cdots (m-n+1)}{n!}$ is the binomial coefficient for $m\in\Z$ and $n\in\Zpos$ and the notation $\NO\cdot\NO$ represents the normal ordering product which is the operation on $\mathcal{A}$ defined by $\NO \psi(n)\NO=\psi(n)$ and 
\begin{align*}
\NO\psi^1(n_1)\cdots\psi^r(n_r)1\NO=
\begin{cases}
\psi^1(n_1)\NO\psi^2(n_2)\cdots\psi^r(n_r)1\NO\quad&\text{if }n_1<0,\\
(-1)^{r-1}\NO\psi^2(n_2)\cdots\psi^r(n_r)1\NO\psi^1(n_1)\quad&\text{if }n_1\geq 0,
\end{cases}
\end{align*}
inductively for $r\in\Zplus$, $n,n_i\in\Z$ and $\psi,\psi^i\in\h$. 
Thus we have a linear map $Y(\,\cdot\,,z):\trip\rightarrow\Hom(\trip,\trip((z)))$ such that $Y(a,z)=\sum_{n\in\Z}a_{(n)}z^{-n-1}$. 
In particular, we have 
\begin{align*}
Y(\psi,z)&=\sum_{n\in\Z}\psi(n)z^{-n-1},\\
Y(a,z)&=\NO\partial^{(n_1-1)}Y(\psi^1,z)\cdots\partial^{(n_r-1)}Y(\psi^r,z)\NO 
\end{align*}
for $a=\psi^1(-n_1)\cdots\psi^r(-n_r)1$, $\psi,\psi^i\in\h$ and $n_i\in\Zplus$, where $\partial^{(n)}=\frac{1}{n!}\frac{d^n}{dz^n}$ for $n\in\Zpos$. 
By definition, $(\psi(-1)1)_{(n)}=\psi(n)$ for any $\psi\in \h$ and $n\in\Z$. 
We may identify $\h$ as a subspace of $\trip$ by the injective map $\psi\mapsto\psi(-1)1+\mathcal{A}_{\geq0}$ and write $\psi_{(m)}$ for $(\psi(-1)1)_{(m)}$ for $\psi\in\h$ and $m\in\Z$. 

Now we set $\1=1+\mathcal{A}_{\geq0}$. 
Since $\psi_{(i)}\1=0$ for $\psi\in\h$ and $i\in\Zpos$, we have $a_{(-1)}\1=a$ and $a_{(i)}\1=0$ for any $a\in\trip$ and $i\in\Zpos$. 
Let $\{(e^i,f^i)\}_{1\leq i\leq d}$ be a canonical basis of $\h$ and consider the vector $\w=\sum_{j=1}^de^j_{(-1)}f^j$. 
We can show that $\w$ does not depend on a choice of a canonical basis. 
One can also see that 
\[
\w_{(0)}\w=\w_{(-2)}\1,\quad\w_{(2)}\w=0,\quad\w_{(1)}\w=2\w\quad\text{and }\w_{(3)}\w=-d\1.
\]
These facts imply that $L_{n}:=\w_{(n+1)}$ with $n\in\Z$ gives a representation of the Virasoro algebra of central charge $-2d$.  
It is easy to check that 
\[
L_0\psi=\psi\quad\text{and } L_i\psi=0\quad\text{for }i\geq 1 
\]
for $\psi\in \h$. 
Hence $\psi\in\h$ is primary of weight $1$ and $[L_n,\psi_{(m)}]=-m\psi_{(m+n)}$ for any $m,n\in\Z$. 
Following \cite{FLM}, we can show the following theorem.  
\begin{theorem}\label{zzol}
The $\mathcal{A}$-module $\trip$ becomes a simple vertex operator superalgebra of central charge $-2d$ with the vacuum vector $\1$ and the Virasoro vector $\w$. 
Furthermore $\trip_0=\C1$ and $\trip=\bigoplus_{n=0}^\infty\trip_n$. 
\end{theorem}

Recall that associated to the $\Z_2$-grading of $\trip$, we have an automorphism $\theta$ of the vertex operator superalgebra $\trip$ of order $2$.   
We denote by $\tripeven$ (resp. $\tripodd$) the $1$-eigenspace (resp. $-1$-eigenspace) for $\theta$. 
Then we see that $\tripeven$ becomes a vertex operator algebra of central charge $-2d$ and $\tripodd$ is an $\tripeven$-module. 
In fact by applying the same arguments in \cite{DM} to a vertex operator superalgebra, we can show the following proposition.  
\begin{proposition}\label{laisaql}
The vertex operator algebra $\tripeven$ is simple, and the $\tripeven$-module $\tripodd$ is irreducible.  
\end{proposition}
We note that $\tripeven=\bigoplus_{n=0}^{\infty}\tripeven_{n}$, $\tripeven_0=\C\1$ and $\tripeven_1=0$. 

From the result in \cite{Lisym} and \cite{Xu}, there exists a unique nondegenerate invariant bilinear form on $\trip$ up to a scalar multiple. 
It is given as follows: 
We consider a bilinear form $(\,\cdot\,,\cdot\,)$ defined by 
\[
(\1,\1)=1,
\]
and
\begin{align*}
(\psi^1_{(-n_1)}\cdots\psi^r_{(-n_r)}\1,\xi^1_{(-m_1)}\cdots\xi^s_{(-m_s)}\1)=\delta_{r,s}(-1)^{\frac{r(r+1)}{2}}\det\left( n_i\delta_{n_i,m_j}\langle\psi^i,\xi^j\rangle\right)_{1\leq i,j\leq r}
\end{align*}
for positive integers $n_i$, $m_j$ and $\psi^i,\xi^i\in\h$. 
We note that $\tripeven$ and $\tripodd$ are mutually orthogonal.
We see that the bilinear form is supersymmetric and hence that the restriction of the bilinear form to $\tripeven$ and $\tripodd$ are symmetric and skew-symmetric respectively. 
We also find that 
\begin{align*}
&(\psi^1_{(-n_1)}\cdots\psi^r_{(-n_r)}\1,\xi^1_{(-m_1)}\cdots\xi^s_{(-m_s)}\1)\\
&=\delta_{r,s}(-1)^{\frac{r(r+1)}{2}}\det\left( n_i\delta_{n_i,m_j}\langle\psi^i,\xi^j\rangle\right)_{1\leq i,j\leq r}\\
&=\delta_{r,s}(-1)^{\frac{r(r+1)}{2}}\sum_{\sigma\in\mathfrak{S}_r}\sgn(\sigma) \prod_{i=1}^r n_i\delta_{n_i,m_{\sigma(i)}}\langle\psi^i,\xi^{\sigma(i)}\rangle\\
&=\sum_{i=1}^s(-1)^{r+i-1}n_1\delta_{n_1,m_i}\langle\psi^1,\xi^i\rangle (\psi^2_{(-n_2)}\cdots\psi^r_{(-n_r)}\1,\xi^1_{(-m_1)}\cdots\xi^{i-1}_{(-m_{i-1})}\xi^{i+1}_{(-m_{i+1})}\cdots \xi^s_{(-m_s)}\1)\\
&=(-1)^{r}(\psi^2_{(-n_2)}\cdots\psi^r_{(-n_r)}\1,\psi^1_{(n_1)}\xi^1_{(-m_1)}\cdots\xi^s_{(-m_s)}\1).
\end{align*}
This proves that $(\psi_{(n)}u,v)=\mp(u,\psi_{(-n)}v)$, that is
\begin{align}\label{uiww;}
(Y(\psi,z)u,v)=\pm (u,Y(e^{zL_1}(-z^{-2})^{L_0}\psi,z^{-1})v)
\end{align}
for any $\psi\in\h$, $u\in\trip^{\pm},v\in\trip$ and $n\in\Z$ respectively. 
Now we have the following proposition. 
\begin{proposition}\label{uurs@pa}
The bilinear form $(\,\cdot\,,\cdot\,)$ on $\trip$ is nondegenerate and invariant. 
\end{proposition} 
\begin{proof}
First we show that the bilinear form $(\,\cdot\,,\cdot\,)$ is nondegenerate. 
Let $\{g^1,\ldots,g^{2d}\}$ and $\{h^1,\ldots,h^{2d}\}$ be bases of $\h$ such that $\langle g^i,h^j\rangle=\delta_{i,j}$ for $1\leq i,j\leq 2d$.
Let $u$ and $v$ be monomials of the form $g^{i_1}_{(-k_1)}\cdots g^{i_t}_{(-k_t)}\1$ and $h^{j_1}_{(-l_1)}\cdots h^{j_r}_{(-l_r)}\1$ with $2\leq i_p,j_p\leq 2d$ and $k_p,l_p\in\Zplus$ respectively. 
Then noting that $g^1_{(n)}v=0$ and $h^1_{(n)}u=0$ for any $n\in\Zpos$, we have 
\[
(g^{1}_{(-n_1)}\cdots g^{1}_{(-n_r)}u,h^{1}_{(-m_1)}\cdots h^{1}_{(-m_s)}v)=\alpha\delta_{r,s}\prod_{i=1}^rn_i\delta_{m_i,n_i}(u,v)
\]
for any $n_i,m_i\in\Zplus$ with $n_1>n_2>\ldots>n_r$, $m_1>m_2>\ldots>m_s$ and some $\alpha\in\{\pm1\}$. 
This implies that for any vector $w$, the coefficient of the monomial $g^{i_1}_{(-k_1)}\cdots g^{i_t}_{(-k_t)}\1$ in $w$ with $1\leq i_j\leq 2d$ and $k_j\in\Zplus$ is a nonzero multiple of the pairing $(w,h^{i_1}_{(-k_1)}\cdots h^{i_t}_{(-k_t)}\1)$. 
Hence if $w$ is in the radical of $(\,\cdot\,,\cdot\,)$ then every coefficient in $w$ of all monomials are zero. 
Thus $w=0$, and this shows that $(\,\cdot\,,\cdot\,)$ is nondegenerate. 

Since $\trip$ is strongly generated by $\h$, \eqref{uiww;} and Proposition \ref{iwiaslka} show that the bilinear form $(\,\cdot\,,\cdot\,)$ is invariant. 
\end{proof}

\subsection{A set of generators of the vertex operator algebra $\tripeven$}\label{Sectb2}
In this section we show that $\tripeven$ is strongly generated by $\tripeven_2\oplus \tripeven_3$. 

We prepare a notation 
\[
B_{m,n}(\psi,\phi)=\frac{(m-1)!(n-1)!}{(m+n-1)!}\psi_{(-m)}\phi_{(-n)}\1
\]
for $\psi,\phi\in\h$ and $m,n\in\Zplus$. 
It is clear that $B_{m,n}(\psi,\phi)=-B_{n,m}(\phi,\psi)$. 
We have the following lemma which will be used frequently. 
\begin{lemma}\label{twpdcnt}
For any $\psi,\phi\in\h$ and $m,n\in\Zplus$ with $m+n\geq 3$, 
\[
B_{m,n}(\psi,\phi)\equiv (-1)^{n-1}B_{m+n-1,1}(\psi,\phi)
\]
modulo the subspace $\haru{L_{-1}B_{m+n-1-i,i}(\psi,\phi)}{1\leq i\leq m+n-2}$. 
\end{lemma}
\begin{proof}
If $n=1$ then there is nothing to prove. 
Suppose that $n\geq1$. 
Then for any $\psi,\phi\in\h$ and $m\in\Zplus$, we have 
\begin{align*}
L_{-1}B_{m,n}(\psi,\phi)&=\frac{(m-1)!(n-1)!}{(m+n-1)!}L_{-1}\psi_{(-m)}\phi_{(-n)}\1\\
&=\frac{m!(n-1)!}{(m+n-1)!}\psi_{(-m-1)}\phi_{(-n)}\1+\frac{(m-1)!n!}{(m+n-1)!}\psi_{(-m)}\phi_{(-n-1)}\1\\
&=(m+n)(B_{m+1,n}(\psi,\phi)+B_{m,n+1}(\psi,\phi)).
\end{align*}
Thus induction on $n$ proves the lemma. 
\end{proof}

We see that $\tripeven_2\oplus \tripeven_3$ is spanned by $B_{m,n}(\psi,\phi)$ for $m,n\in\Zplus$ with $m+n=2,3$ and $\psi,\phi\in\h$.
Set 
\[
U:=\haru{a^1_{(-n_1)}\cdots a^s_{(-n_s)}\1}{a^i\in \tripeven_2\oplus \tripeven_3,\,n_i\in\Zplus}.
\] 
The aim of this section is to prove that $\tripeven=U$. 
To show this we first prove the following lemma. 
\begin{lemma}\label{oiai}
For any $m,n\in\Zplus$ and $\psi,\phi\in\h$, $B_{m,n}(\psi,\phi)\in U$.
\end{lemma}
\begin{proof} 
Let $\{(e^i,f^i)\}_{1\leq i\leq d}$ be a canonical basis of $\h$.
It suffices to show Lemma \ref{oiai} for the pairs $(\psi,\phi)=(e^i,f^j)$, $(e^i,e^j)$ and $(f^i,f^j)$ with $1\leq i,j\leq d$. 
We first prove that the vectors $B_{m,n}(e^i,f^j)$ and $B_{m,n}(e^i,e^j)$ lie in $U$ for any $m,n\in\Zplus$ and $1\leq i,j\leq d$.
To show this we use induction on $m+n$. 
It is clear that the lemma holds if $m+n\leq 3$. 
Let $m,n\in\Zplus$ with $m+n\geq 4$ and assume that $B_{p,q}(e^i,f^j), B_{p,q}(e^i,e^j)\in U$ for any $p,q\in\Zplus$ with $p+q<m+n$ and $1\leq i,j\leq d$. 
Since $L_{-1}U\subset U$, induction hypothesis and Lemma \ref{twpdcnt} show that 
\begin{align}\label{pper}
B_{m,n}(\psi,\phi)\equiv (-1)^{n-1}B_{m+n-1,1}(\psi,\phi)\mod U. 
\end{align}

For simplicity, we set $\psi^j=e^j$ or $f^j$. 
Then by using \eqref{lqqljq}, we calculate  
\begin{align*}
B_{1,1}(e^i,f^i)_{(-1)}B_{1,q}(e^{i},\psi^j)&=(q+2)\delta_{i,j}B_{1,q+2}(e^i,\psi^j)+\binom{q+2}{2}B_{3,q}(e^i,\psi^j) 
\end{align*}
for any $q\in\Zplus$.
Thus by \eqref{pper} and induction hypothesis we have  
\[
\binom{m+n-1+\delta_{i,j}}{2}B_{3,m+n-3}(e^i,\psi^j)=B_{1,1}(e^i,f^i)_{(-1)}B_{1,m+n-3}(e^{i},\psi^j)\in U.
\]
This proves that $B_{3,m+n-3}(e^i,\psi^j)\in U$, and hence $B_{m,n}(e^i,\psi^j)\in U$ by \eqref{pper}. 

By exchanging the canonical basis $\{(e^i,f^i)\}_{1\leq i\leq d}$ to $\{(f^i,-e^i)\}_{1\leq i\leq d}$, we find that $B_{m,n}(f^i,f^j)\in U$ for any $1\leq i,j\leq d$. 
\end{proof}

Now we consider the subspaces defined by ${\filt}^0\trip=\C\1$ and 
\[
{\filt}^r\trip=\haru{\psi^1_{(-n_1)}\cdots\psi^s_{(-n_s)}\1}{\psi^i\in\h,n_i\in\Zplus\text{ and }s\leq r}.
\]
for $r\in\Zplus$. 
Then we have a filtration $\{{\filt}^r\trip\}_{r\in\Zpos}$ on $\trip$.
We set 
\[
{\filt}^r\tripeven={\filt}^r\trip\cap\tripeven
\]
for $r\in\Zpos$.   
It is easy to see that  
\begin{align*}
{\filt}^{2r+1}\tripeven=&{\filt}^{2r}\tripeven,\\
{\filt}^r\tripeven=&\bigoplus_{k=0}^\infty {\filt}^r\tripeven_k\quad\text{with }{\filt}^r\tripeven_k:={\filt}^r\tripeven\cap\tripeven_k
\end{align*}
for $r\in\Zpos$. 

We note that for any $j,r\in\Zpos$ and $\psi\in\h$, 
\begin{align}\label{aikiadbca}
\psi_{(j)}{\filt}^r\trip\subset {\filt}^{r-1}\trip.
\end{align}
This implies that 
\begin{align}\label{iwikwdc}
a_{(j)}{\filt}^r\trip\subset {\filt}^{r}\trip\quad\text{for any $j,r\in\Zpos$  and $a\in {\filt}^2\tripeven$}.
\end{align}
By using \eqref{iwikwdc}, we can show the following lemma. 
\begin{lemma}\label{qqse}
Fix $r\in\Zpos$ and suppose that ${\filt}^r\tripeven\subset U$. 
Then for any $m\in\Zplus$ and $a\in{\filt}^2\tripeven$, $a_{(-m)}{\filt}^r\tripeven\subset U$. 
\end{lemma}
\begin{proof}
We set
\[
W_k=\haru{a_{(-m)}u}{a\in {\filt}^2\tripeven_s\text{ with }s\leq k, u\in {\filt}^r\tripeven,m\in\Zplus}
\] 
for $k\geq 2$. 
We shall show that $W_k\subset U$ for any $k\geq 2$ by using induction on $k$.
Since ${\filt}^r\tripeven\subset U$ by the assumption, we see that $W_3\subset U$.
Let $k>3$ and assume that $W_s\subset U$ for any $s<k$.
We notice from the proof of Lemma \ref{oiai} that for any $a\in {\filt}^2\tripeven_k$ with $k\geq4$ there exist $u^i\in\tripeven_2$, $v^i\in{\filt}^2\tripeven_{k-2}$ and $w\in{\filt}^2\tripeven_{k-1}$ such that 
\begin{align*}
a=\sum_{i}u^i_{(-1)}v^i+L_{-1}w. 
\end{align*}
For any $m\in\Zplus$ and $u\in {\filt}^r\tripeven$, the identity \eqref{supertwBorcherds} with $p=-1$, $q=0$ and $r=-m$ leads 
\begin{align*}
a_{(-m)}u=mw_{(-m-1)}u+\sum_{i}\sum_{j=0}^\infty(u^i_{(-1-j)}v^i_{(-m+j)}u+v^i_{(-m-1-j)}u^i_{(j)}u).    
\end{align*}
Thus by \eqref{iwikwdc} and induction hypothesis, one sees that the vectors of the form  $w_{(-m-1)}u$, $v^i_{(-m+j)}u$ and $v^i_{(-m-1-j)}u^i_{(j)}u$ are in $U$ for any $i$ and $j\in\Zpos$. 
Since $u^i\in \tripeven_2$ we get $a_{(-m)}u\in U$. 
This shows $W_k\subset U$.  
\end{proof}

Now we can show the following proposition. 
\begin{proposition}\label{llsp}
The vertex operator algebra $\tripeven$ is strongly generated by $\tripeven_2\oplus\tripeven_3$. 
\end{proposition}
\begin{proof}
It is enough to show that for any $r\in\Zpos$, ${\filt}^r\tripeven\subset U$. 
We use induction on $r$. 
We see that ${\filt}^0\tripeven=\C\1\subset U$ and ${\filt}^2\tripeven\subset U$ by Lemma \ref{oiai}.
Let $r\geq 2$ and assume that ${\filt}^{2r-2}\tripeven\subset U$. 
By \eqref{lqqljq} and \eqref{aikiadbca}, we see that  
\begin{align}\label{ujuhjsdja}
\psi^1_{(-n_1)}\cdots\psi^{2r}_{(-n_{2r})}\1\equiv (\psi^1_{(-n_1)}\psi^2_{(-n_2)}\1)_{(-1)}\psi^3_{(-n_3)}\cdots\psi^{2r}_{(-n_{2r})}\1
\end{align}
modulo ${\filt}^{2r-2}\tripeven$ for any $r\in\Zpos$, $\psi^i\in\h$ and $n_i\in\Zplus$. 
Since $\psi^3_{(-n_3)}\cdots\psi^{2r}_{(-n_{2r})}\1\in {\filt}^{2r-2}\tripeven$, Lemma \ref{qqse} and induction hypothesis prove that the right hand side in \eqref{ujuhjsdja} is in $U$. 
Thus so is the left hand side in \eqref{ujuhjsdja}. 
This implies that ${\filt}^{2r}\tripeven\subset U$. 
\end{proof}

We set 
\begin{align*}
e^{i,j}&:=B_{1,1}(e^i,e^j)=e^{i}_{(-1)}e^j,\\
h^{i,j}&:=B_{1,1}(e^i,f^j)=e^{i}_{(-1)}f^j,\\
f^{i,j}&:=B_{1,1}(f^i,f^j)=f^{i}_{(-1)}f^j,\\
E^{i,j}&:=B_{2,1}(e^i,e^j)+B_{2,1}(e^j,e^i)=\frac{1}{2}(e^i_{(-2)}e^j+e^j_{(-2)}e^i),\\
H^{i,j}&:=B_{2,1}(e^i,f^j)+B_{2,1}(f^j,e^i)=\frac{1}{2}(e^i_{(-2)}f^j+f^j_{(-2)}e^i),\\
F^{i,j}&:=B_{2,1}(f^i,f^j)+B_{2,1}(f^j,f^i)=\frac{1}{2}(f^i_{(-2)}f^j+f^j_{(-2)}f^i)
\end{align*}
for any $1\leq i,j \leq d$. 
We note that $e^{i,j},h^{i,j},f^{i,j}$ are quasi-primary and $E^{i,j},H^{i,j},F^{i,j}$ are primary. 
Since $L_{-1}\psi_{(-1)}\phi=\psi_{(-2)}\phi-\phi_{(-2)}\psi$ for $\psi,\phi\in \h$,
we see that 
\begin{align*}
\tripeven_2&=\bigoplus_{i,j=1}^d \C h^{i,j}\oplus\bigoplus_{1\leq i<j \leq d}(\C e^{i,j}\oplus\C f^{i,j}),\\
\tripeven_3&=\bigoplus_{i,j=1}^d\C H^{i,j}\oplus\bigoplus_{1\leq i\leq j \leq d}(\C E^{i,j}\oplus\C F^{i,j})\oplus L_{-1}\tripeven_2.
\end{align*}
Thus we have the following corollary.
\begin{corollary}\label{alaiwdcql}
Let $\{(e^i,f^i)\}_{1\leq i\leq d}$ be a canonical basis. 
Then $\tripeven$ is strongly generated by the vectors $e^{i,j}$, $h^{i,j}$, $f^{i,j}$, $E^{i,j}$, $H^{i,j}$ and $F^{i,j}$ with $1\leq i,j \leq d$. 
\end{corollary}
In the case $d=1$, we see that $\tripeven$ is strongly generated by $h^{1,1}=\w$, $E:=E^{1,1}$, $H:=H^{1,1}$ and $F:=F^{1,1}$. 
It is proved in \cite{K1} that the OPEs among $\w,E,H,F$ are coincides with that of the triplet algebra with $c=-2$.  

\subsection{$C_2$-cofiniteness of the vertex operator algebra $\tripeven$}\label{Sectb3}
In this section we shall show that the vertex operator algebra $\tripeven$ is $C_2$-cofinite. 

We first consider the case $d=1$. 
We denote by $\T$ the vertex operator superalgebra $\trip$ with $d=1$,
and set $\T^\pm=\trip^\pm$ respectively.  
We give a proof of the following theorem. 
\begin{theorem}\label{aiugyqwo}
The vertex operator algebra $\Teven$ is $C_2$-cofinite. 
\end{theorem} 
This theorem follows from Proposition \ref{alalia} below. 
\begin{remark}
However it is not asserted in \cite{GaKa} that the vacuum representation of $c=-2$ triplet algebra is $C_2$-cofinite, we can show this fact by means of the explicit forms of null vectors although the calculations for null vectors and relations among generators are not so easy.  
We here give a direct proof of $C_2$-cofiniteness of $\Teven$ and find explicit relations among generators of the commutative algebra $\Teven/C_2(\Teven)$.  
\end{remark}

We now start proving the following proposition.  
\begin{proposition}\label{alalia}
The dimension of $\Teven/C_2(\Teven)$ is less than or equal to $11$.
\end{proposition} 

Let $\{(e,f)\}$ be a canonical basis of $\h$. 
By Corollary \ref{alaiwdcql}, $\Teven$ is strongly generated by 
\[
\w=e_{(-1)}f,\quad E:=e_{(-2)}e,\quad H:=\frac{1}{2}(e_{(-2)}f+f_{(-2)}e)\quad\text{and }F:=f_{(-2)}f. 
\] 
Thus $\Teven/C_2(\Teven)$ is generated as an algebra by $\overline{w},\overline{E},\overline{H}$ and $\overline{F}$ by Proposition \ref{iwiwsoi}. 
We shall find relations among these generators. 
For simplicity, we reset 
\[
A_{m,n}(\psi,\phi)=(m+n-1)!B_{m,n}(\psi,\phi)=(m-1)!(n-1)!\psi_{(-m)}\phi_{(-n)}\1
\]
for $\psi,\phi\in\h$ and $m,n\in\Zplus$. 
Then by the proof of Lemma \ref{twpdcnt}, we have  
\begin{align*}
A_{m,n}(\psi,\phi) \equiv (-1)^{n-1}A_{m+n-1,1}(\psi,\phi) \mod L_{-1}\Teven. 
\end{align*}
Since $L_{-1}\Teven\subset C_2(\Teven)$, we see that 
\begin{align}\label{iaclia}
\overline{A_{m,n}(\psi,\phi)}= (-1)^{n-1}\overline{A_{m+n-1,1}(\psi,\phi)} 
\end{align}
for any $\psi,\phi\in\h$ and $m,n\in\Zplus$. 
We set 
\[
\Gamma_m(\psi,\phi)=\overline{A_{m-1,1}(\psi,\phi)}\in \Teven/C_2(\Teven)
\]
for any $\psi,\phi\in\h$ and $m\geq 2$. 
We note that
\[
\Gamma_m(\psi,\phi)=\overline{A_{m-1,1}(\psi,\phi)}=-\overline{A_{1,m-1}(\phi,\psi)}=(-1)^{m-1}\Gamma_m(\phi,\psi)
\]
for any $\psi,\phi\in\h$ and $m\geq 2$. 
In particular, $\Gamma_{2m}(\psi,\psi)=0$ for any positive integer $m\in\Zplus$. 

To find relations in $\Teven/C_2(\Teven)$, we use the identity  
\begin{align}
\begin{split}\label{laisjcqa}
\Gamma_m(\psi,\phi)\cdot\Gamma_k(\xi,\eta)&=(m-2)!(k-2)!\overline{\psi_{(-m+1)}\phi_{(-1)}\xi_{(-k+1)}\eta}\\
&\quad+\frac{\langle\phi,\xi\rangle}{k}\Gamma_{m+k}(\psi,\eta)-(-1)^k\frac{\langle\phi,\eta\rangle}{2}\Gamma_{m+k}(\psi,\xi)\\
&\quad+(-1)^{m+k}\frac{\langle\psi,\eta\rangle}{m}\Gamma_{m+k}(\phi,\xi)-(-1)^{m}\frac{\langle\psi,\xi\rangle}{m+k-2}\Gamma_{m+k}(\phi,\eta) 
\end{split}
\end{align}
for $\psi,\phi,\xi,\eta\in\h$ and $m,k\geq 2$, which can be proved by a direct calculation using the associativity formula and \eqref{iaclia}. 

First we have  
\begin{align}\label{ipde}
\Gamma_3(e,e)\cdot\Gamma_3(e,e)=\Gamma_3(f,f)\cdot\Gamma_3(f,f)=0.
\end{align}
Next \eqref{laisjcqa} gives the relation 
\begin{align*}
\Gamma_3(f,e)\cdot\Gamma_4(e,e)=-\frac{2}{15}\Gamma_{7}(e,e). 
\end{align*} 
Since $\Gamma_4(e,e)=0$, we see that $\Gamma_{7}(e,e)=0$. 
Therefore, 
\begin{align}\label{ppsr}
\Gamma_2(f,e)^2\cdot\Gamma_3(e,e)=\frac{7}{12}\Gamma_7(e,e)=0. 
\end{align}
As well, we get 
\begin{align}\label{ppsr3}
\Gamma_2(e,f)^2\cdot\Gamma_3(f,f)=\frac{7}{12}\Gamma_7(f,f)=0. 
\end{align}
If we take a canonical basis $\{(\frac{1}{\sqrt{2}}(e+f),\frac{1}{\sqrt{2}}(e-f))\}$, then we have 
\[
\Gamma_2(e,f)^2\cdot\Gamma_3\left(\frac{1}{\sqrt{2}}(e+f),\frac{1}{\sqrt{2}}(e+f)\right)=0
\]
because $\Gamma_2(e,f)=-\Gamma_2\left(\frac{1}{\sqrt{2}}(e+f),\frac{1}{\sqrt{2}}(e-f)\right)$. 
Since 
\[
\Gamma_3\left(\frac{1}{\sqrt{2}}(e+f),\frac{1}{\sqrt{2}}(e+f)\right)=\frac{1}{2}\left(\Gamma_3(e,e)+2\Gamma_3(e,f)+\Gamma_3(f,f)\right),
\]
one has 
\begin{align}\label{ppsr2}
\Gamma_2(e,f)^2\cdot \Gamma_3(e,f)=0. 
\end{align}

Next we use the identities  
\begin{align*}
\Gamma_m(e,e)\cdot\Gamma_k(f,f)&=(m-2)!(k-2)!\overline{e_{(-m+1)}e_{(-1)}f_{(-k+1)}f}\\
&\quad+\left(-\frac{1}{k}+\frac{(-1)^{k}}{2}+\frac{(-1)^m}{m+k-2}-\frac{(-1)^{m+k}}{m}\right)\Gamma_{m+k}(e,f)
\end{align*}
and 
\begin{align*}
&\Gamma_m(e,f)\cdot\Gamma_k(f,e)\\
&=(m-2)!(k-2)!\overline{e_{(-m+1)}f_{(-1)}f_{(-k+1)}e}-\frac{(-1)^{k}}{2}\Gamma_{m+k}(e,f)+\frac{(-1)^m}{m+k-2}\Gamma_{m+k}(f,e)\\
&=-(m-2)!(k-2)!\overline{e_{(-m+1)}e_{(-1)}f_{(-k+1)}f}-\left(\frac{(-1)^{k}}{2}+\frac{(-1)^k}{m+k-2}\right)\Gamma_{m+k}(e,f)
\end{align*} 
for $m,k\geq 2$, where we note $\Gamma_{m+k}(f,e)=(-1)^{m+k-1}\Gamma_{m+k}(e,f)$.

Since $(-1)^{k-1}\Gamma_k(f,f)=\Gamma_k(f,f)$ for $k\geq 2$, we get the following relation:
\begin{align}\begin{split}\label{alsio}
&\Gamma_m(e,f)\cdot\Gamma_k(e,f)\\
&=-\Gamma_m(e,e)\cdot\Gamma_k(f,f)+\left(\frac{(-1)^{m}}{m}+\frac{(-1)^k}{k}+\frac{1-(-1)^{m+k}}{m+k-2}\right)\Gamma_{m+k}(e,f).
\end{split}\end{align}
On the other hand \eqref{laisjcqa} shows that 
\begin{align}\label{liqevlq}
\Gamma_m(e,f)\cdot\Gamma_k(e,f)=\left(\frac{1}{m}+\frac{1}{k}\right)\Gamma_{m+k}(e,f)
\end{align}
for $k,m\geq 2$. 
One now has 
\begin{align*}
&\Gamma_3(e,f)^2=-\Gamma_3(e,e)\cdot\Gamma_3(f,f)-\frac{2}{3}\Gamma_6(e,f),\quad \Gamma_2(e,f)^3=\frac{3}{4}\Gamma_6(e,f)
\end{align*}
by \eqref{alsio}. 
It also follows from \eqref{liqevlq} that 
\begin{align}\label{hhea}
\Gamma_3(e,f)^2=\frac{2}{3}\Gamma_6(e,f)=\frac{8}{9}\Gamma_2(e,f)^3.
\end{align}
Therefore, we have 
\begin{align}\label{oosnr}
\Gamma_3(e,e)\cdot\Gamma_3(f,f)=-\frac{16}{9}\Gamma_2(e,f)^3.
\end{align}

We here recall that
\begin{align*}
\overline\w&=\overline{B_{1,1}(e,f)}=\Gamma_2(e,f),\\
\overline{E}&=(\overline{B_{2,1}(e,e)}+\overline{B_{2,1}(e,e)})=\Gamma_3(e,e),\\
\overline{H}&=(\overline{B_{2,1}(e,f)}+\overline{B_{2,1}(f,e)})=\Gamma_3(e,f),\\
\overline{F}&=(\overline{B_{2,1}(e,f)}+\overline{B_{2,1}(f,e)})=\Gamma_3(f,f).
\end{align*}
Then from \eqref{ipde}--\eqref{ppsr2}, \eqref{hhea} and \eqref{oosnr}, we have the following relations:
\begin{align}
&\overline{E}\cdot\overline{E}=\overline{F}\cdot\overline{F}=\overline{H}\cdot \overline{E}=\overline{H}\cdot \overline{F}=0,\label{aadoc1}\\
&2\overline{H}^2=-\overline{E}\cdot \overline{F}=\frac{16}{9}\overline{\w}^3,\label{aadoc2}\\
&\overline{w}^2\cdot \overline{E}=\overline{w}^2\cdot \overline{H}=\overline{w}^2\cdot \overline{F}=0.\label{aadoc3}
\end{align}
Now we can prove Proposition \ref{alalia}. 
\vskip1ex
\noindent
\textit{ Proof of Proposition \ref{alalia}.} 
We recall that $\Teven/C_2(\Teven)$ is generated by $\overline\w,\overline{E},\overline{H}$ and $\overline{F}$ as a commutative algebra. 
Note that \eqref{aadoc2} and \eqref{aadoc3} give $\overline{w}^5=0$. 
Hence by \eqref{aadoc1}--\eqref{aadoc3}, we see that $\Teven/C_2(\Teven)$ is spanned by $11$ vectors 
\[
\overline{\w}^i\quad\text{for }i=0,1,2,3,4,\quad\overline{E},\quad\overline{\w}\cdot\overline{E},\quad\overline{H},\quad\overline{\w}\cdot\overline{H},\quad\overline{F},\quad\overline{\w}\cdot\overline{F}.
\] 
Hence $\dim \Teven/C_2(\Teven)\leq 11$. 
\hfill$\Box$
\vskip1ex

By using Theorem \ref{aiugyqwo} and Proposition \ref{rrlsiw} we have 
\begin{theorem}\label{laisuqwi}
The vertex operator algebra $\tripeven$ is $C_2$-cofinite. 
\end{theorem} 
\begin{proof}
Let $\{(e^i,f^i)\}_{1\leq i\leq d}$ be a canonical basis  of $\h$.
Then $\h=\bigoplus_{i=1}^{d}(\C e^i\oplus\C f^i)$ is a direct sum of mutually orthogonal subspaces. 
This orthogonal decomposition induces an embedding of a vertex operator algebra $\bigotimes_{i=1}^{d}\Teven$ in $\tripeven$.
By the definition of the Virasoro vector of $\tripeven$, it is clear that $\tripeven$ contains a full vertex operator subalgebra isomorphic to $\bigotimes_{i=1}^{d}\Teven$.
On the other hand, it is known that a tensor product of $C_2$-cofinite vertex operator algebras is also $C_2$-cofinite.
Therefore, $\bigotimes_{i=1}^{d}\Teven$ is $C_2$-cofinite by Theorem \ref{aiugyqwo}. 
Then Proposition \ref{rrlsiw} shows that $\trip^+$ is $C_2$-cofinite.
\end{proof}

\section{Classification of irreducible $\tripeven$-modules}\label{Sectc}
In this section we classify irreducible $\tripeven$-modules.  
Firstly we construct an irreducible $\theta$-twisted $\trip$-module. 
Secondly we show, by using Zhu's algebra, that any irreducible $\tripeven$-module appears in the irreducible $\trip$-module $\trip$ or in the irreducible $\theta$-twisted $\trip$-module. 
\subsection{Irreducible $\tripeven$-modules}\label{Sectc1}
We have shown that $\trip^{\pm}$ are irreducible $\tripeven$-modules. 
To find the other irreducible $\tripeven$-modules we construct a $\theta$-twisted $\trip$-module. 
The construction can be done as in the case of the free bosonic vertex operator algebra. 
Set 
\[
\hat{L}^{\theta}(\h):=\h\otimes t^{\frac{1}{2}}\C[t^{\pm1}]\oplus\C K
\]
and make it a superspace such that the even part is $\C K$ and the odd part is $\h\otimes t^{\frac{1}{2}}\C[t^{\pm1}]$. 
Then $\hat{L}^{\theta}(\h)$ has a Lie superalgebra structure by 
\[
[\psi\otimes t^m,\phi\otimes t^n]_+=m\delta_{m+n,0}\langle\psi,\phi\rangle K,\quad[K,\hat{L}^{\theta}(\h)]=0
\]
for $\psi,\phi\in\h$ and $m,n\in\frac{1}{2}+\Z$. 
Now we consider the associative algebra $\mathalpha{A}^{\theta}$ which is the quotient algebra of the universal enveloping algebra $\mathcal{U}(\hat{L}^\theta(\h))$ by the ideal generated by $K-1$. 
Canonically, the algebra $\mathalpha{A}^{\theta}$ has a $\Z_2$-grading $\mathalpha{A}^{\theta}=\mathalpha{A}^{\theta}_{\bar{0}}\oplus\mathalpha{A}^{\theta}_{\bar{1}}$ of algebras.
Thus $\mathalpha{A}^{\theta}$ has naturally an involution which is also denoted by $\theta$.  

Let $\mathcal{A}^\theta_{>0}$ be the left ideal of $\mathalpha{A}^{\theta}$ generated by the vectors $\psi\otimes t^m $ for $\psi\in\h$ and $m\in\frac{1}{2}+\Zpos$, and set $\triptw=\mathcal{A}^\theta/\mathcal{A}^\theta_{>0}$. 
We set $1_{\theta}=1+\mathcal{A}^\theta_{>0}\in \triptw$ and denote by $\psi{(m)}$ the operator of left multiplication by $\psi\otimes t^m$ on $\triptw$ for $\psi\in\h$ and $m\in\frac{1}{2}+\Z$. 
Since the involution $\theta$ preserves $\mathcal{A}^\theta_{>0}$, the $\Z_2$-grading of $\mathcal{A}^\theta$ induces a decomposition $\triptw=\triptweven\oplus\triptwodd$, where $\triptw^\pm$ are $\pm1$-eigenspaces for $\theta$ of $\triptw$. 
We note that $\triptw\cong\Lambda(\h\otimes t^{-\frac{1}{2}}\C[t^{-1}])$ as vector spaces. 

We can endowed $\triptw$ with a $\theta$-twisted $\trip$-module structure following \cite{FLM}. 
First we set 
\[
W(\psi,z)=\sum_{i\in\frac{1}{2}+\Z}\psi(i)z^{-i-1}
\]
for any $\psi\in\h$ and define 
\[
W(v,z)=\NO \partial^{(n_1-1)}W(\psi^1,z)\cdots \partial^{(n_r-1)}W(\psi^r,z)\NO
\]
for $v=\psi^1_{(-n_1)}\cdots\psi^r_{(-n_r)}\1$ for $\psi^i\in\h$ and $n_i\in\Zplus$. 
This defines a linear map $W(\,\cdot\,,z)$ from $\trip$ to $\Hom(\triptw,\triptw((z^{\frac{1}{2}})))$. 
We take the coefficients $c_{m,n}\in\C\,(m,n\in\Zpos)$ subject to the formal expansion
\[
\sum_{m,n\geq0} c_{mn}x^{m}y^{n}=-\log \left(
\frac{(1+x)^{\frac{1}{2}}+(1+y)^{\frac{1}{2}}}{2}\right),
\]
and consider the operator 
\begin{align*}
\Delta(z)&=\sum_{m,n\geq 0}\sum_{i=1}^{d} c_{mn}
{e^{i}}_{(m)} {f^{i}}_{(n)}z^{-m-n}-\sum_{m,n\geq 0}\sum_{i=1}^{d} c_{mn}
 {f^{i}}_{(m)}{e^{i}}_{(n)}z^{-m-n}\\
&=2\sum_{m,n\geq 0}\sum_{i=1}^{d} c_{mn}
{e^{i}}_{(n)}{f^{i}}_{(m)} z^{-m-n} 
\end{align*}
on $\trip$, where $\{(e^i,f^i)\}_{1\leq i\leq d}$ is a canonical basis of $\h$.
The last equality holds because $c_{m,n}=c_{n,m}$ for any $m,n\in \Zpos$. 
By using the operator $\Delta(z)$, the vertex operator associated to the vector $v=\psi^1_{(-n_1)}\cdots\psi^r_{(-n_r)}\1$ is defined by 
\[
Y(v,z)=W(e^{\Delta(z)}v,z),
\]
and the pair $(\triptw,Y(\,\cdot\,,z))$ becomes a $\theta$-twisted $\trip$-module. 
Since $\Omega(\triptw)=\C1_{\theta}$, $\triptw$ is an irreducible $\trip$-module. 
We see that $\triptw^\pm$ become $\tripeven$-modules. 
By applying \cite[Theorem 5.5]{DLi} to a vertex operator superalgebra, we have 
\begin{proposition}\label{alihala}
The $\tripeven$-modules $\triptw^{\pm}$ are irreducible as $\tripeven$-modules.
\end{proposition} 

We write $\psi_{(m)}$ for $(\psi_{(-1)}\1)_{(m)}$ for any $\psi\in\h$ and $m\in\frac{1}{2}+\Z$. 
Direct calculations show that 
\begin{align*} 
Y(\psi,z)&=W(\psi,z),\\
Y(\psi_{(-1)}\phi,z)&=W(\psi_{(-1)}\phi,z)+\frac{\langle\psi,\phi\rangle}{8}\id\,z^{-2}
\end{align*}
for $\psi,\phi\in\h$. 
Since $[L_{m},\psi_{(n)}]=-m\psi_{(m+n)}$ for $m\in\Z$ and $n\in\frac{1}{2}+\Z$, we see that 
\begin{align*}
L_{0}1_{\theta}&=-\frac{d}{8}1_{\theta},\\
L_{0}\psi^1_{(-n_1)}\cdots\psi^r_{(-n_r)}1_{\theta}&=\left(-\frac{d}{8}+\sum n_i\right)\psi^1_{(-n_1)}\cdots\psi^r_{(-n_r)}1_{\theta}
\end{align*} 
for any $\psi^i\in\h$ and $n_i\in\frac{1}{2}+\Zpos$. 
Therefore we find that 
\[
\triptw=\bigoplus_{i=0}^\infty\triptw_{-\frac{d}{8}+\frac{i}{2}}.
\]
In fact we have 
\[
\triptweven=\bigoplus_{i=0}^\infty\triptw_{-\frac{d}{8}+i},\quad \triptwodd=\bigoplus_{i=0}^\infty\triptw_{\frac{-d+4}{8}+i}.
\]

We here state how does $o(a)=a_{(1)}$ act on $\Omega(M)$ for $a\in\tripeven_2$ and the  known irreducible $\tripeven$-modules $M$.
We fix a canonical basis $\{(e^i,f^i)\}_{1\leq i\leq d}$. 
Then $\tripeven_2$ is spanned by the vectors $h^{i,j}$, $e^{i,j}$ and $f^{i,j}$ for $1\leq i,j\leq d$. 
It is easy to see that $o(h^{i,j})=o(e^{i,j})=o(f^{i,j})=0$ on $\Omega(\tripeven)=\C\1$ and that $o(h^{i,j})=-\frac{1}{8}\delta_{i,j}\id$ and $o(e^{i,j})=o(f^{i,j})=0$ on $\Omega(\triptweven)=\C 1_\theta$. 
The spaces $\Omega(\tripodd)$ and $\Omega(\triptwodd)$ have basis $\{e^i,f^i\}_{1\leq i\leq d}$ and $\{e^i_{(-\frac{1}{2})}1_\theta, f^i_{(-\frac{1}{2})}1_\theta\}_{1\leq i\leq d}$ respectively.
We set $x_1^i=e^i$, $y_1^i=f^i$, $x_{\frac{3}{8}}^i=e^i_{(-\frac{1}{2})}1_\theta$ and $y_{\frac{3}{8}}^i=f^i_{(-\frac{1}{2})}1_\theta$ for $1\leq i\leq d$. 
Then we can calculate that 
\begin{alignat}{4}
o(h^{i,i})x_{h}^k&=h\delta_{i,k}x_h^{i},&\qquad o(h^{i,i})y_{h}^k&=h\delta_{i,k}y_h^{i},\label{eedort1}\\
o(h^{i,j})x_{h}^k&={\alpha_h}\delta_{j,k}x_h^{i},&\qquad o(h^{i,j})y_{h}^k&=-{\alpha_h}\delta_{i,k}y_h^{j},\label{eedort2}\\
o(e^{i,j})x_{h}^k&=0,&\qquad o(e^{i,j})y_{h}^k&=-{\alpha_h}(\delta_{j,k}x_h^{i}+\delta_{i,k}x_h^{j}),\label{eedort3}\\
o(f^{i,j})x_{h}^k&={\alpha_h}(\delta_{j,k}y_h^{i}+\delta_{i,k}x_h^{j}),&\qquad o(f^{i,j})y_{h}^k&=0\label{eedort4}
\end{alignat}
for $h=1,\frac{3}{8}$ and $1\leq i, j,k\leq d$ with $i\neq j$, where we set $\alpha_1=1$ and $\alpha_{\frac{3}{8}}=\frac{1}{2}$.  
As we will show later, Zhu's algebra $A(\tripeven)$ is generated by the images of vectors of weight $2$ when $d\geq 2$ (see Proposition \ref{qloacql}). 
Thus \eqref{eedort1}--\eqref{eedort4} characterize the known irreducible $A(\tripeven)$-modules in the case $d\geq 2$.     

\subsection{Main Theorem}\label{Main}
In this section we state the main theorem, and describe some products of vectors in Zhu's algebra $A(\tripeven)$. 

The main theorem in the paper is  
\begin{theorem}\label{ikwksldl}
The list $\{\trip^{\pm},\triptw^\pm\}$ gives a complete list of inequivalent irreducible $\tripeven$-modules. 
\end{theorem}
By Theorem \ref{kiuwewytt}, to prove Theorem \ref{ikwksldl}, we only have to classify irreducible $A(\tripeven)$-modules. 
The classification of irreducible modules of $A(\tripeven)$ will be given in Section \ref{Sectc2} and Section \ref{Sectc3}, where we treat the cases $d=1$ and $d\geq2$ respectively. 

Before considering the structure of $A(\tripeven)$ in these two cases, we shall prepare notations and give some identities in $A(\tripeven)$. 
As in the proof of Lemma \ref{twpdcnt}, for $\psi,\phi\in\h$ and $m,n\in\Zplus$, $(m+n)\left(B_{m+1,n}(\psi,\phi)+B_{m,n+1}(\psi,\phi)\right)=L_{-1}B_{m,n}(\psi,\phi)$. 
Since $[L_{-1}a]=-[L_0a]$ in $A(\tripeven)$ for any $a\in \tripeven$, we see that  
\begin{align}\label{aosjncqa}
[B_{m+1,n}(\psi,\phi)]+[B_{m,n+1}(\psi,\phi)]=-[B_{m,n}(\psi,\phi)]. 
\end{align}  
We set 
\[
\Theta_{m}(\psi,\phi):=[B_{m-1,1}(\psi,\phi)]=\frac{1}{m-1}[\psi_{(-m+1)}\phi]\in A(\tripeven)
\] 
for $\psi,\phi\in\h$ and $m\geq 2$.
Then we have the following lemma.
\begin{lemma}\label{alidjqal}
For any $\psi,\phi\in\h$ and $m,n\in\Zplus$,
\[
[B_{m,n}(\psi,\phi)]=(-1)^{n-1}\sum_{j=0}^{n-1}\binom{n-1}{j}\Theta_{m+n-j}(\psi,\phi).
\]
In particular, for any $\psi,\phi\in\h$ and integer $m\geq 2$,
\[
\Theta_m(\phi,\psi)=(-1)^{m-1}\sum_{i=0}^{m-2}\binom{m-2}{i}\Theta_{m-i}(\psi,\phi).
\] 
\end{lemma}
\begin{proof}
We use induction on $p=m+n$ and \eqref{aosjncqa}.
The case $p=2$ is clear. 
Let $p\geq3 $ and suppose that the lemma is valid when $m+n<p$. 
We prove that the lemma holds for $m=p-i$ and $n=i$ for $1\leq i\leq p-1$ by using induction on $i$. 
By definition, $[B_{p-1,1}(\psi,\phi)]=\Theta_{p}(\psi,\phi)$. 
This shows the lemma for $m=p-1,n=1$.
Let $i\geq 1$. 
By using \eqref{aosjncqa}, we have $[B_{p-i-1,i+1}(\psi,\phi)]=-[B_{p-i,i}(\psi,\phi)]-[B_{p-i-1,i}(\psi,\phi)]$.  
Thus  
\begin{align*}
[B_{p-i-1,i+1}(\psi,\phi)]=&(-1)^{i}\sum_{j=0}^{i-1}\binom{i-1}{j}\Theta_{p-j}(\psi,\phi)+(-1)^{i}\sum_{j=0}^{i-1}\binom{i-1}{j}\Theta_{p-1-j}(\psi,\phi)\\
=&(-1)^{i}\sum_{j=0}^{i}\binom{i}{j}\Theta_{p-j}(\psi,\phi). 
\end{align*}
This proves the first identity in the lemma for $m=p-i-1$ and $n=i+1$. 
The second identity follows from the first identity because $\Theta_m(\phi,\psi)=-[B(\psi,\phi)_{1,m-1}]$.  
\end{proof}

We state some identities which can be shown by direct calculations and by using Lemma \ref{alidjqal}. 
We will use these identities in the next two sections.
For $\psi,\phi,\xi,\eta\in\h$ and $m\geq 2$, 
\begin{align}
\begin{split}\label{aliuhcq}
&\Theta_{2}(\psi,\phi)*\Theta_{m}(\xi,\eta)\\
&=\frac{1}{m-1}[\psi_{(-1)}\phi_{(-1)}\xi_{(-m+1)}\eta]\\
&\quad+\langle\phi,\xi\rangle((m+1)\Theta_{m+2}(\psi,\eta)+2m\Theta_{m+1}(\psi,\eta)+(m-1)\Theta_{m}(\psi,\eta))\\
&\quad-\langle\psi,\xi\rangle((m+1)\Theta_{m+2}(\phi,\eta)+2m\Theta_{m+1}(\phi,\eta)+(m-1)\Theta_{m}(\phi,\eta))\\
&\quad+\langle\phi,\eta\rangle\left(\binom{m+1}{2}\Theta_{m+2}(\xi,\psi)+2\binom{m}{2}\Theta_{m+1}(\xi,\psi)+\binom{m-1}{2}\Theta_{m}(\xi,\psi)\right)\\
&\quad-\langle\psi,\eta\rangle\left(\binom{m+1}{2}\Theta_{m+2}(\xi,\phi)+2\binom{m}{2}\Theta_{m+1}(\xi,\phi)+\binom{m-1}{2}\Theta_{m}(\xi,\phi)\right),
\end{split}\end{align}
and  
\begin{align}
\begin{split}\label{aliuhalsiuuq}
\Theta_{3}(\psi,\phi)*\Theta_{m}(\xi,\eta)
&=\frac{1}{2(m-1)}[\psi_{(-2)}\phi_{(-1)}\xi_{(-m+1)}\eta]\\
&\quad+\langle\phi,\xi\rangle\left( \binom{m+2}{2}\Theta_{m+3}(\psi,\eta)+3\binom{m+1}{2}\Theta_{m+2}(\psi,\eta)\right.\\
&\left.\qquad+3\binom{m}{2}\Theta_{m+1}(\psi,\eta)+\binom{m-1}{2}\Theta_{m}(\psi,\eta)\right)\\
&\quad+\frac{m\langle\psi,\xi\rangle}{2}((m+2)\Theta_{m+3}(\phi,\eta)+3(m+1)\Theta_{m+2}(\phi,\eta)\\
&\qquad+3m\Theta_{m+1}(\phi,\eta)+(m-1)\Theta_{m}(\phi,\eta))\\
&\quad-\frac{3\langle\phi,\eta\rangle}{2}\left(\binom{m+2}{3}\Theta_{m+3}(\xi,\psi)+m\binom{m+1}{2}\Theta_{m+2}(\xi,\psi)\right.\\
&\qquad\left.+m\binom{m}{2}\Theta_{m+1}(\xi,\psi)+\binom{m}{3}\Theta_{m}(\xi,\psi)\right)\\
&\quad-\langle\psi,\eta\rangle\left(\binom{m+2}{3}\Theta_{m+3}(\xi,\phi)+3\binom{m+1}{3}\Theta_{m+2}(\xi,\phi)\right.\\
&\qquad\left.+3\binom{m}{3}\Theta_{m+1}(\xi,\phi)+\binom{m-1}{3}\Theta_{m}(\xi,\phi)\right).
\end{split}\end{align}
 
\subsection{Classification of irreducible $A(\tripeven)$-modules for $d=1$}\label{Sectc2}
We keep using the notation $\Teven$ for the vertex operator algebra $\tripeven$ with $d=1$. 
In this section we classify irreducible $A(\Teven)$-modules. 
To do this we will find enough relations in Zhu's algebra of $\Teven$ (see \cite{GaKa}). 

We take a canonical basis $\{(e,f)\}$ of $\h$. 
We first remark that $\Theta_2(e,f)=[\w]$ is in the center of $A(\Teven)$.  
It is also clear that 
\begin{align}
\Theta_3(e,e)*\Theta_m(e,e)&=0,\quad\Theta_3(f,f)*\Theta_m(f,f)=0\label{rrks}
\end{align}
for any $m\geq 2$.  

Lemma \ref{alidjqal} implies that 
\begin{align}
\Theta_2(e,e)&=0\label{laijsdsq},\\
\Theta_4(e,e)&=-\Theta_3(e,e),\label{liuhcwq}\\
\begin{split}
\Theta_6(e,e)&=-2\Theta_5(e,e)-3\Theta_4(e,e)-2\Theta_3(e,e)=-2\Theta_5(e,e)+\Theta_3(e,e).\label{aoscql;qa}
\end{split}
\end{align}
We apply \eqref{aliuhalsiuuq} for $m=3$ and $4$ and use \eqref{laijsdsq}--\eqref{aoscql;qa} to get  
\begin{align}
\Theta_3(f,e)*\Theta_3(e,e)&=11\Theta_5(e,e)-10\Theta_3(e,e),\label{iiwmj}\\  
\Theta_3(f,e)*\Theta_4(e,e)&=-8\Theta_{7}(e,e)+12\Theta_5(e,e)-5\Theta_{3}(e,e).
\end{align}
Since $\Theta_3(f,e)*(\Theta_3(e,e)+\Theta_4(e,e))=0$ by \eqref{liuhcwq}, one has 
\begin{align}\label{soiuw;laq}
\Theta_7(e,e)=\frac{23}{8}\Theta_5(e,e)-\frac{15}{8}\Theta_3(e,e).
\end{align}
Moreover, \eqref{aliuhcq}, \eqref{laijsdsq}--\eqref{soiuw;laq} give   
\begin{align}
\Theta_2(e,f)*\Theta_3(e,e)&=10\Theta_5(e,e)-9\Theta_3(e,e),\label{alishcqa}\\
\Theta_2(e,f)*\Theta_5(e,e)&=\frac{83}{8}\Theta_5(e,e)-\frac{75}{8}\Theta_3(e,e).
\end{align}
Hence we find that 
\begin{align*}
\begin{split}
\Theta_2(e,f)^2*\Theta_3(e,e)&=10\Theta_2(e,f)*\Theta_5(e,e)-9\Theta_2(e,f)*\Theta_3(e,e)\\
&=\frac{55}{4}\Theta_5(e,e)-\frac{51}{4}\Theta_3(e,e).
\end{split}
\end{align*}
Since 
\begin{align}
\Theta_5(e,e)=\frac{1}{10}\Theta_2(e,f)*\Theta_3(e,e)+\frac{9}{10}\Theta_3(e,e)
\end{align}
by \eqref{alishcqa} we have the identity 
\begin{align}
(\Theta_2(e,f)-1)*(8\Theta_2(e,f)-3)*\Theta_3(e,e)=0.\label{lalwjncl} 
\end{align}
Moreover \eqref{iiwmj} shows that 
\[
\Theta_3(f,e)*\Theta_3(e,e)=\frac{11}{10}\Theta_2(e,f)*\Theta_3(e,e)-\frac{1}{10}\Theta_3(e,e).
\]
Therefore noting
\begin{align}
\Theta_3(e,f)=\Theta_3(f,e)-\Theta_2(e,f)
\label{ooer}
\end{align} 
which follows from Lemma \ref{alidjqal}, we have 
\begin{align}\label{ttwods}
(\Theta_3(e,f)+\Theta_3(f,e))*\Theta_3(e,e)=\frac{1}{5}(6\Theta_2(e,f)-1)*\Theta_3(e,e). 
\end{align}

By exchanging the choice of the canonical basis $\{(e,f)\}$ to $\{(f,-e)\}$, we get the relations   
\begin{align}
&(\Theta_2(e,f)-1)*\left(\Theta_2(e,f)-\frac{3}{8}\right)*\Theta_3(f,f)=0,\label{aqwqwoil}\\
&(\Theta_3(e,f)+\Theta_3(f,e))*\Theta_3(f,f)=-\frac{1}{5}(6\Theta_2(e,f)-1)*\Theta_3(f,f).\label{uunek} 
\end{align}
If we exchange the canonical basis $\{(e,f)\}$ to $\{(\frac{1}{\sqrt{2}}(e+f),\frac{1}{\sqrt{2}}(e-f))\}$ and use the relations \eqref{lalwjncl}--\eqref{uunek} with respect to the latter basis, then we get 
\begin{align}
&(\Theta_2(e,f)-1)*\left(\Theta_2(e,f)-\frac{3}{8}\right)*(\Theta_3(e,f)+\Theta_3(f,e))=0.\label{ccede}
\end{align} 

Next we calculate 
\begin{align}
\Theta_2(e,f)^2&=6\Theta_{4}(e,f)+6\Theta_{3}(e,f)+\Theta_{2}(e,f),\label{aliuhq1}\\
\Theta_2(e,f)*\Theta_3(e,f)&=10\Theta_{5}(e,f)+12\Theta_{4}(e,f)+3\Theta_{3}(e,f),\label{aouwhcqpw}\\
\Theta_2(e,f)*\Theta_4(e,f)&=15\Theta_6(e,f)+20\Theta_{5}(e,f)+6\Theta_4(e,f).\label{sdpivadn}
\end{align}
by using \eqref{aliuhcq}. 
Hence $\Theta_4(e,f)$ and $\Theta_5(e,f)$ are expressed by using $\Theta_2(e,f)$ and $\Theta_3(e,f)$ as 
\begin{align}
\Theta_4(e,f)&=-\Theta_3(e,f)+\frac{1}{6}\Theta_2(e,f)^2-\frac{1}{6}\Theta_2(e,f),\label{oasqlcq}\\
\Theta_5(e,f)&=\frac{1}{10}\Theta_2(e,f)*\Theta_3(e,f)+\frac{9}{10}\Theta_3(e,f)-\frac{1}{5}\Theta_2(e,f)^2+\frac{1}{5}\Theta_2(e,f). \label{8rweqw9}
\end{align}
Substituting \eqref{oasqlcq} and \eqref{8rweqw9} into the both sides in \eqref{sdpivadn} yields 
\begin{align}
\begin{split}\label{laislaisdb}
\Theta_6(e,f)=&-\frac{1}{5}\Theta_2(e,f)*\Theta_3(e,f)-\frac{4}{5}\Theta_3(e,f)\\
&+\frac{1}{90}\Theta_2(e,f)^3+\frac{17}{90}\Theta_2(e,f)^2-\frac{1}{5}\Theta_2(e,f).
\end{split}
\end{align}

Now we calculate that   
\begin{align}
\begin{split}\label{laIqaow0}
\Theta_3(f,f)*\Theta_3(e,e)=&-\frac{1}{4}[f_{(-2)}e_{(-1)}e_{(-2)}f]\\
&-\frac{85}{2}\Theta_6(e,f)-73\Theta_5(e,f)-36\Theta_4(e,f)-\frac{9}{2}\Theta_3(e,f),
\end{split}\\
\begin{split}\label{laIqaow2}
\Theta_3(f,e)*\Theta_3(e,f)=&\frac{1}{4}[f_{(-2)}e_{(-1)}e_{(-2)}f]\\
&+\frac{45}{2}\Theta_6(e,f)+45\Theta_5(e,f)+27\Theta_4(e,f)+\frac{9}{2}\Theta_3(e,f),
\end{split}\\
\begin{split}\label{laIqaow3}
\Theta_3(e,f)^2=&20\Theta_6(e,f)+30\Theta_5(e,f)+12\Theta_4(e,f)+\Theta_3(e,f)
\end{split}
\end{align} 
by using \eqref{aliuhalsiuuq} for $m=3$ and Lemma \ref{alidjqal}.
By \eqref{ooer} we have  
\begin{align*}
(\Theta_3(e,f)+\Theta_3(f,e))^2&=(2\Theta_3(e,f)+\Theta_2(e,f))^2\\
&=4\Theta_3(e,f)^2+\Theta_2(e,f)*(4\Theta_3(e,f)+\Theta_2(e,f)).
\end{align*}
Hence by \eqref{laIqaow3} and \eqref{oasqlcq}--\eqref{laislaisdb} we can show the relation   
\begin{align}\label{uuwosw}
(\Theta_3(e,f)+\Theta_3(f,e))^2=\frac{1}{9}\Theta_2(e,f)^2*(8\Theta_2(e,f)+1).
\end{align}
We also find that by using \eqref{ooer}, 
\begin{align}
\begin{split}\label{ccie}
(\Theta_3(e,f)+\Theta_3(f,e))^2&=(2\Theta_3(f,e)-\Theta_2(e,f))*(2\Theta_3(e,f)+\Theta_2(e,f))\\
&=4\Theta_3(f,e)*\Theta_3(e,f)+\Theta_2(e,f)^2
\end{split}
\end{align}
by \eqref{ooer}. 
On the other hand, by \eqref{oasqlcq}--\eqref{laIqaow2}, we see that 
\begin{align*}
&4\Theta_3(f,f)*\Theta_3(e,e)+4\Theta_3(f,e)*\Theta_3(e,f)+\Theta_2(e,f)^2\\
&=-80\Theta_6(e,f)-112\Theta_5(e,f)-36\Theta_4(e,f)+\Theta_2(e,f)^2\\
&=\frac{2}{5}(6\Theta_2(e,f)-1)*(2\Theta_3(e,f)+\Theta_2(e,f))-\frac{1}{9}\Theta_2(e,f)^2*(8\Theta_2(e,f)+1). 
\end{align*}
Therefore, by \eqref{uuwosw} and \eqref{ccie}, we have 
\begin{align}
\begin{split}\label{wwurn}
&4\Theta_3(f,f)*\Theta_3(e,e)\\
&=\frac{2}{5}(6\Theta_2(e,f)-1)*(\Theta_3(e,f)+\Theta_3(f,e))-\frac{2}{9}\Theta_2(e,f)^2*(8\Theta_2(e,f)+1).
\end{split} 
\end{align}
Finally \eqref{ccede} and \eqref{uuwosw} imply the identity  
\begin{align}\label{urwls}
\Theta_2(e,f)^2*(8\Theta_2(e,f)+1)*(\Theta_2(e,f)-1)*(8\Theta_2(e,f)-3)=0.
\end{align}

By definition we see that 
\[
[\w]=\Theta_2(e,f),\quad [E]=2\Theta_3(e,e),\quad [H]=\Theta_3(e,f)+\Theta_3(f,e)\quad\text{and }[F]=2\Theta_3(f,f). 
\]
Since $\w$, $E$, $H$ and $F$ are quasi-primary, we have $\Phi([\w])=[\w]$ and $\Phi([x])=-[x]$ for $x=E,H,F$.  
Thus by the identities \eqref{rrks}, \eqref{lalwjncl}, \eqref{ttwods}--\eqref{ccede}, \eqref{uuwosw}--\eqref{urwls} and by applying the anti-involution $\Phi$ to these identities, we have the following proposition: 
\begin{proposition}\label{yypwgtd}
The vector $[\w]$ is a central element and satisfies the relations  
\begin{align*}
[\w]^2*(8[\w]+1)*([\w]-1)*(8[\w]-3)&=0,\\
([\w]-1)*(8[\w]-3)*[x]&=0
\end{align*}
for $x=E,H$ and $F$.
The relations  
\begin{align*}
[E]^2&=[F]^2=0,\quad [H]^2=\frac{1}{9}[\w]^2*(8[\w]+1),\\
[H]*[E]&=-[E]*[H]=\frac{1}{5}(6[\w]-1)*[E],\\
[H]*[F]&=-[F]*[H]=-\frac{1}{5}(6[\w]-1)*[F],\\
[E]*[F]&=-\frac{2}{5}(6[\w]-1)*[H]-\frac{2}{9}[\w]^2*(8[\w]+1),\\
[F]*[E]&=\frac{2}{5}(6[\w]-1)*[H]-\frac{2}{9}[\w]^2*(8[\w]+1)
\end{align*}
also hold in $A(\Teven)$. 
\end{proposition}
\begin{remark}
We can show that $\g=\C E\oplus \C H\oplus\C F=\{x\in \Teven_3\,|\,L_1x=0\}$ becomes a Lie algebra with commutation relation $[x,y]=\frac{1}{5}p(x_{(2)}y)$ for $x,y\in \g$, where $p$ denotes the projection from $\Teven$ to the subspace of all quasi-primary vectors. 
In fact, $\g$ is isomorphic to $\spl$ such that $[E,F]=-2H,\,[H,E]=E$ and $[H,F]=-F$ ($\{iE,2H,iF\}$ forms a standard basis). 
It is easy to see that the bilinear form $(\,\cdot\,,\cdot)$ gives a nondegenerate invariant symmetric bilinear form of $\g$. 
We note that $(H,H)=1$, $(E,F)=-2$.
Therefore the bilinear form coincides with the twice of the normalized Killing form. 
We can then unify the relations in Proposition \ref{yypwgtd} with respect to $E,H,F$ to the relations
\begin{align*}
&([\w]-1)*(8[\w]+3)*[x]=0,\\
&[x]*[y]=\frac{1}{5}(6[\w]-1)*[[x,y]]+\frac{(x,y)}{9}[\w]^2*(8[\w]+1)
\end{align*}
for $x,y\in\g$. 
These imply that any $A(\Teven)$-module is a module for $\spl$ (see \cite{GaKa}).   
\end{remark}
It follows from Proposition \ref{eert} and Corollary \ref{alaiwdcql} that Zhu's algebra $A(\Teven)$ is generated by $[\w]$, $[E]$, $[H]$ and $[F]$. 
By using Proposition \ref{yypwgtd}, we shall find some ideals of $A(\Teven)$.  
We now set 
\begin{align}
v_{0}:&=-\frac{1}{9}\left(13[\w]-3\right)*\left([\w]-1\right)*\left(8[\w]+1\right)*\left(8[\w]-3\right),\\
w_0:&=[\w]*\left([\w]-1\right)*\left(8[\w]+1\right)*\left(8[\w]-3\right),\\
v_{-\frac{1}{8}}:&=\frac{128}{9}[\w]^2*\left([\w]-1\right)*\left(8[\w]-3\right).
\end{align}
Then we have 
\begin{align}
v_{k}*v_{l}&=\delta_{k,l}v_{k}\quad\text{for }k,l=0,-1/8,\\
v_{0}*w_0&=w_0*v_{0}=w_0,\quad w_0^2=0,\\
v_{-\frac{1}{8}}*w_0{}&=w_0*v_{-\frac{1}{8}}=0.
\end{align}
We note that $[\w]^2*v_1=[\w]*w_0=0$ and $[\w]*v_{-\frac{1}{8}}=-\frac{1}{8}v_{-\frac{1}{8}}$. 
It is also valid that $v_k*[x]=w_0*[x]=0$ for $x=E,H,F$. 
Therefore $A_0=\C v_0+\C w_0$ and $A_{-\frac{1}{8}}=\C v_{-\frac{1}{8}}$ are ideals of $A(\Teven)$.   

We now consider the subspace $A_1$ spanned by the vectors $[\w]^2*\left(8[\w]+1\right)*\left(8[\w]-3\right)$ and $\left(8[\w]-3\right)*[x]$ with $x=E,H$ and $F$.  
We also denote by $A_{\frac{3}{8}}$ the subspace spanned by the vectors $[\w]^2*\left(8[\w]+1\right)*\left([\w]-1\right)$ and $\left([\w]-1\right)*[x]$ with $x=E,H$ and $F$.
It is clear from Proposition \ref{yypwgtd} that $[w]$ acts on $A_\lambda$ by the scalar $\lambda$ for $\lambda=1,\frac{3}{8}$. 
Proposition \ref{yypwgtd} also implies that 
\begin{align}\label{sdfrt}
A(\Teven)=A_0\oplus A_1\oplus A_{-\frac{1}{8}}\oplus A_{\frac{3}{8}}. 
\end{align}
and that $A_{\lambda}*A_{\mu}=0$ if $\lambda\neq \mu$ for $\lambda,\mu=0,1,-\frac{1}{8},\frac{3}{8}$. 
Therefore, the decomposition \eqref{sdfrt} is a direct sum of ideals. 
\begin{proposition}\label{laisdqli}
The ideals $A_\lambda$ for $\lambda=1,\frac{3}{8}$ are homomorphic images of the $2\times 2$ matrix algebra $M_2(\C)$. 
\end{proposition}
\begin{proof}
We set 
\begin{alignat*}{4}
A=&\frac{1}{90}[\w]^2*\left(8[\w]+1\right)*\left(8[\w]-3\right),&\quad B&=\frac{1}{10}(8[\w]-3)*[H],\\
C=&\frac{1}{10}(8[\w]-3)*[E],&\quad D&=\frac{1}{10}(8[\w]-3)*[F]
\end{alignat*}
when $\lambda=1$ and set 
\begin{alignat*}{4}
A=&-\frac{64}{45}[\w]^2*\left(8[\w]+1\right)*\left([\w]-1\right),&\quad B&=-\frac{16}{5}([\w]-1)*[H],\\
C=&-\frac{16}{5}([\w]-1)*[E],&\quad D&=-\frac{16}{5}([\w]-1)*[F].
\end{alignat*}
when $\lambda=\frac{3}{8}$. 
Then by Proposition \ref{yypwgtd}, we have Table \ref{iieewa}, which is the multiplication tables among $A,B,C$ and $D$.  

\begin{table}[htbp]
\caption{The table of $x*y$ for $x,y=A,B,C,D$}
\begin{center}
\begin{tabular}{c|cccc}  \label{iieewa}  
$x\backslash y$&$A$&$B$&$C$&$D$\\
\hline
$A$&$\frac{1}{2}A$&$\frac{1}{2}B$&$\frac{1}{2}C$&$\frac{1}{2}D$\\   
$B$&$\frac{1}{2}B$&$\frac{1}{2}A$& $\frac{1}{2}C$&$-\frac{1}{2}D$\\
$C$&$\frac{1}{2}C$&$-\frac{1}{2}C$&$0$&$-A-B$\\
$D$&$\frac{1}{2}D$&$\frac{1}{2}D$&$-A+B$&$0$\\
\end{tabular}
\end{center}
\end{table}

Hence for $\lambda=1,\frac{3}{8}$, if we set  
\begin{align*}
v_\lambda^{1,1}&=-A-B,\quad v_\lambda^{1,2}=C,\quad v_\lambda^{2,1}=D,\quad v_\lambda^{2,2}=-A+B,  
\end{align*} 
then we get the identities $v_\lambda^{i,j}*v_\lambda^{k,l}=\delta_{j,k}v_\lambda^{i,l}$ for any $i,j,k,l=1,2$.  
This shows that $A_\lambda$ is a homomorphic image of $M_2(\C)$ for $\lambda=1,\frac{3}{8}$. 
\end{proof}
\begin{remark}\label{remark4}
Let $W$ be an $A(\Teven)$-module and $\lambda=1$ or $\frac{3}{8}$. 
We recall the scalars $\alpha_1=1$ and $\alpha_{\frac{3}{8}}=\frac{1}{2}$ defined in the last paragraph in Section \ref{Sectc1}. 
For any $w\in W$ such that $[\w].w=\lambda w$, we have 
\begin{alignat*}{4}
v_\lambda^{1,1}.w&=-\frac{1}{2}\left(w+\frac{1}{\alpha_\lambda^2}[H].w\right),&\quad v_\lambda^{1,2}.w&=\frac{1}{2\alpha_\lambda^2}[E].w,\\
 v_\lambda^{2,1}.w&=\frac{1}{2\alpha_\lambda^2}[F].w,&\quad v_\lambda^{2,2}.w&=-\frac{1}{2}\left(w-\frac{1}{\alpha_\lambda^2}[H].w\right). 
\end{alignat*} 
Therefore, the nonzero vectors $x_1:=v_\lambda^{1,1}.w$ and $x_2:=v_\lambda^{2,1}.w$ satisfy $v_\lambda^{i,j}.x_k=\delta_{j,k}x_i$ if there exists $0\neq w\in W$ with $[\w].w=\lambda w$ for $\lambda=1,\frac{3}{8}$. 
\end{remark}
As a corollary of Proposition \ref{laisdqli} we can show the following theorem.
\begin{theorem}\label{uhealasc}
Zhu's algebra $A(\Teven)$ has only four inequivalent irreducible modules $\Omega(\T^\pm)$ and $\Omega(\Ttw^\pm)$. 
\end{theorem}  
\begin{proof}
Let $W$ be an irreducible $A(\Teven)$-module. 
Then $W=A_\lambda.w$ for some $\lambda=0,1,-\frac{1}{8},\frac{3}{8}$ and a nonzero vector. 
We note that $A_0$ is commutative and the other ideals are homomorphic image of simple algebras. 
Thus there exist at most four irreducible $A(\Teven)$ module. 
On the other hand, there exist four inequivalent irreducible modules $\Omega(\T^\pm)$ and $\Omega(\Ttw^\pm)$. 
Therefore we have the theorem.  
\end{proof}
It is still possible that the ideal $A_0$ degenerates to one dimensional. 
The fact that $A_0$ is just two dimensional can be proved by showing the existence of a reducible indecomposable $A(\Teven)$-module on which $[\w]$ does not act diagonally but nilpotently. 
We will show that such $A(\Teven)$-modules indeed exist (see Remark \ref{wwes} below). 
Then we see that $\dim A(\Teven)\geq 11$. 
Since $\dim A(\Teven)\leq \dim \Teven/C_2(\Teven)$ as in Section \ref{Secta2}, by Proposition \ref{alalia}, we have 
\[
\dim A(\Teven)=\dim \Teven/C_2(\Teven)=11. 
\]  

\subsection{Classification of irreducible $A(\tripeven)$-modules for $d\geq 2$}\label{Sectc3}
In this section we classify irreducible $A(\tripeven)$-modules in the case $d>1$ (see Theorem \ref{liawqwqq} below).
First we note that  
\begin{align*}
[e^{i,j}]&=\Theta_2(e^i,e^j),\quad [h^{i,j}]=\Theta_2(e^i,f^j),\quad [f^{i,j}]=\Theta_2(f^i,f^j),\\
[E^{i,j}]&=\Theta_3(e^i,e^j)+\Theta_3(e^j,e^i)=2\Theta_3(e^i,e^j)+\Theta_2(e^i,e^j),\\
[H^{i,j}]&=\Theta_3(e^i,f^j)+\Theta_3(f^j,e^i)=2\Theta_3(e^i,f^j)+\Theta_2(e^i,f^j),\\
[F^{i,j}]&=\Theta_3(f^i,f^j)+\Theta_3(f^j,f^i)=2\Theta_3(f^i,f^j)+\Theta_2(f^i,f^j)
\end{align*}
for $1\leq i,j\leq d$. 
By Corollary \ref{alaiwdcql} and Proposition \ref{eert}, we see that Zhu's algebra $A(\tripeven)$ is generated by the vectors above. 
Moreover, we can show the following proposition. 
\begin{proposition}\label{qloacql} If $d>1$, then Zhu's algebra $A(\tripeven)$ is generated by $[e^{i,j}]$, $[h^{i,j}]$ and $[f^{i,j}]$ for $1\leq i,j\leq d$. 
\end{proposition}
\begin{proof}
Let $1\leq i,j\leq d$. 
It suffices to express $[E^{i,j}]$, $[H^{i,j}]$ and $[F^{i,j}]$ by means of $[e^{k,l}]$, $[h^{k,l}]$ and $[f^{k,l}]$ with $1\leq k,l\leq d$. 
We shall assume that $i\neq j$. 
By using \eqref{aliuhcq} and Lemma \ref{alidjqal}, we can calculate that  
\begin{align}
[h^{i,i}]*[e^{i,j}]&=3\Theta_4(e^i,e^j)+4\Theta_3(e^i,e^j)+\Theta_2(e^i,e^j),\label{uiahysdc1}\\
[h^{j,j}]*[e^{i,j}]&=3\Theta_4(e^i,e^j)+2\Theta_3(e^i,e^j),\label{uiahysdc2}\\
[h^{i,i}]*[h^{i,j}]&=3\Theta_4(e^i,f^j)+4\Theta_3(e^i,f^j)+\Theta_2(e^i,f^j),\label{uiahysdc3}\\
[h^{j,j}]*[h^{i,j}]&=3\Theta_4(e^i,f^j)+2\Theta_3(e^i,f^j),\label{uiahysdc4}\\
[h^{i,i}]*[f^{i,j}]&=3\Theta_4(f^i,f^j)+4\Theta_3(f^i,f^j)+\Theta_2(f^i,f^j),\label{uiahysdc5}\\
[h^{j,j}]*[f^{i,j}]&=3\Theta_4(f^i,f^j)+2\Theta_3(f^i,f^j).\label{uiahysdc6}
\end{align}
Hence we see that 
\begin{align}
[E^{i,j}]&=2\Theta_3(e^i,e^j)+\Theta_2(e^i,e^j)=([h^{i,i}]-[h^{j,j}])*[e^{i,j}],\label{lakshjf1}\\
[H^{i,j}]&=2\Theta_3(e^i,f^j)+\Theta_2(e^i,f^j)=([h^{i,i}]-[h^{j,j}])*[h^{i,j}],\\
[F^{i,j}]&=2\Theta_3(f^i,f^j)+\Theta_2(f^i,f^j)=([h^{i,i}]-[h^{j,j}])*[f^{i,j}].
\end{align}

It also follows from \eqref{aliuhcq} and a similar formula with $f^i$ instead of $e^i$ that    
\begin{align*}
[h^{i,j}]*[e^{i,j}]=3\Theta_4(e^i,e^i)+4\Theta_3(e^i,e^i),\quad [h^{j,i}]*[f^{i,j}]=3\Theta_4(f^i,f^i)+2\Theta_3(f^i,f^i).
\end{align*}
Therefore, by \eqref{liuhcwq}, we have 
\begin{align}\label{ualru}
[E^{i,i}]=2\Theta_{3}(e^i,e^i)=2[h^{i,j}]*[e^{i,j}],\quad [F^{i,i}]=2\Theta_{3}(f^i,f^i)=-2[h^{j,i}]*[f^{i,j}].
\end{align}

Finally we calculate that $[e^{i,j}]*[f^{i,j}]$ and $[h^{i,j}]*[h^{j,i}]$ for $1\leq i\neq j\leq d$. 
By \eqref{aliuhcq} and Lemma \ref{alidjqal}, we find 
\begin{align*}
[e^{i,j}]*[f^{i,j}]&=[e^i_{(-1)}e^j_{(-1)}f^i_{(-1)}f^j]+3\Theta_4(e^i,f^i)+4\Theta_3(e^i,f^i)+\Theta_2(e^i,f^i)\\
&\quad+3\Theta_4(e^j,f^j)+4\Theta_3(e^j,f^j)+\Theta_2(e^j,f^j),\\
[h^{i,j}]*[h^{j,i}]&=[e^i_{(-1)}f^j_{(-1)}e^j_{(-1)}f^i]+3\Theta_4(e^i,f^i)+4\Theta_3(e^i,f^i)+\Theta_2(e^i,f^i)\\
&\quad+3\Theta_4(e^j,f^j)+2\Theta_3(e^j,f^j).
\end{align*}
Then \eqref{oasqlcq} and the fact that $\Theta_3(e^k,f^k)=\frac{1}{2}([H^{k,k}]-[h^{k,k}])$ for $1\leq k\leq d$ show that 
\begin{align}
[e^{i,j}]*[f^{i,j}]&=\frac{1}{2}([H^{i,i}]+[H^{j,j}])+\frac{1}{2}([h^{i,i}]-[h^{j,j}])^2,\label{owq;autrf1}\\
[h^{i,j}]*[h^{j,i}]&=\frac{1}{2}([H^{i,i}]-[H^{j,j}])+\frac{1}{2}([h^{i,i}]-[h^{j,j}])^2.\label{owq;autrf2}
\end{align}
These imply that 
\begin{align}
[H^{j,j}]&=[e^{i,j}]*[f^{i,j}]-[h^{i,j}]*[h^{j,i}],\label{jmreeari1}\\
[H^{i,i}]&=[e^{i,j}]*[f^{i,j}]+[h^{i,j}]*[h^{j,i}]-([h^{i,i}]-[h^{j,j}])^2.\label{jmreea,ri2}
\end{align}
This completes the proof.
\end{proof}
Apply the anti-automorphism $\Phi$ to the both sides in \eqref{lakshjf1}--\eqref{owq;autrf1}. 
We note that the anti-involution $\Phi$ acts trivially on $[e^{i,j}]$, $[h^{i,j}]$ and $[f^{i,j}]$ and by the scalar $-1$ on $[E^{i,j}]$, $[H^{i,j}]$ and $[F^{i,j}]$ for $1\leq i,j\leq d$. 
Thus we have  
\begin{align}
[h^{i,i}]*[x^{i,j}]&=[x^{i,j}]*[h^{j,j}]\quad \text{for }x^{i,j}=e^{i,j},\,h^{i,j},\,h^{j,i},\,f^{i,j},\label{pqkelo1}\\
[h^{i,j}]*[y^{i,j}]&=-[y^{i,j}]*[h^{i,j}]\quad \text{for }y^{i,j}=e^{i,j},\,f^{i,j}\label{pqkelo6},\\
[e^{i,j}]*[f^{i,j}]&=-[f^{i,j}]*[e^{i,j}]+([h^{i,i}]-[h^{j,j}])^2 \label{pqkelo7}
\end{align}
for any $1\leq i\neq j\leq d$.  
We also see that \eqref{jmreeari1} and \eqref{jmreea,ri2} prove the relation 
\begin{align}\label{owklwqa}
[h^{i,j}]*[h^{j,i}]+[h^{j,i}]*[h^{i,j}]=([h^{i,i}]-[h^{j,j}])^2
\end{align} 
for $1\leq i\neq j\leq d$. 
 
It is clear that $[h^{i,i}]$ and $[h^{j,j}]$ for $1\leq i,j\leq d$ commute with each other in $A(\tripeven)$. 
Since $A(\tripeven)$ is finite dimensional, irreducible $A(\tripeven)$-modules are of finite dimension. 
Thus any irreducible $\tripeven$-module is a sum of simultaneous generalized eigenspaces for the actions of $[h^{i,i}]$ for all $1\leq i\leq d$. 
We show that the actions of $[h^{i,i}]$ with $1\leq i\leq d$ are diagonal on any irreducible $A(\tripeven)$-module. 
\begin{proposition}\label{utnugurnde}
Let $W$ be an irreducible $A(\tripeven)$-module.
Then $W$ is a direct sum of simultaneous eigenspaces for the actions of $[h^{i,i}]$ with all $1\leq i\leq d$. 
\end{proposition} 
\begin{proof}
Since $W$ is finite dimensional, $W$ contains a nonzero simultaneous eigenvector $w$ for the actions of $[h^{i,i}]\,(1\leq i\leq d)$. 

Consider the vector of the form $[e^{k,l}].w$ for $1\leq k\neq l\leq d$. 
By \eqref{uiahysdc1}, \eqref{uiahysdc2}, \eqref{pqkelo1} and a trivial calculation, one has 
\begin{align*}
[h^{k,k}].[e^{k,l}].w&=[e^{k,l}].[h^{l,l}].w=\lambda_l[e^{k,l}].w,\\
[h^{l,l}].[e^{k,l}].w&=[e^{k,l}].[h^{k,k}].w=\lambda_k[e^{k,l}].w,\\
[h^{i,i}].[e^{k,l}].w&=[e^{k,l}].[h^{i,i}].w=\lambda_i[e^{k,l}].w\quad\text{for }i\neq k, l,
\end{align*} 
where $\lambda_i$ is the eigenvalue of $w$ for the action of $[h^{i,i}]$. 
Hence, we see that $[e^{k,l}].w$ is also simultaneous eigenvector for the actions of $[h^{i,i}]$ with $1\leq i\leq d$. 
As well, \eqref{uiahysdc3}--\eqref{uiahysdc6} prove that $[h^{k,l}].w$ and $[f^{k,l}].w$ are also simultaneous eigenvectors for the actions of $[h^{i,i}]\,(1\leq i\leq d)$. 
Finally Proposition \ref{qloacql} implies that $A(\tripeven).w$ is a direct sum of simultaneous eigenspaces for the actions of $[h^{i,i}]$ for all $1\leq i\leq d$.
Since $W$ is irreducible, we have $W=A(\tripeven).w$.
Hence the proposition holds. 
\end{proof}

We now find more relations in $A(\tripeven)$. 
We recall the relation
\[
[H^{i,i}]*[E^{i,i}]=\frac{1}{5}(6[h^{i,i}]-1)*[E^{i,i}]
\]
for $1\leq i\leq d$ in Proposition \ref{yypwgtd}.
By \eqref{ualru}, \eqref{jmreeari1} and \eqref{pqkelo1}--\eqref{pqkelo7}, for $1\leq j\leq d$ with $j\neq i$,  
\begin{align*}
[H^{i,i}]*[E^{i,i}]&=([e^{i,j}]*[f^{i,j}]-[h^{j,i}]*[h^{i,j}])*(2[h^{i,j}]*[e^{i,j}])\\
&=2[e^{i,j}]*[f^{i,j}]*[h^{i,j}]*[e^{i,j}]\\
&=2[h^{i,j}]*[e^{i,j}]*[f^{i,j}]*[e^{i,j}]\\
&=2([h^{i,i}]-[h^{j,j}])^2*[h^{i,j}]*[e^{i,j}]\\
&=([h^{i,i}]-[h^{j,j}])^2*[E^{i,i}]
\end{align*}
because $[e^{i,j}]^2=[h^{i,j}]^2=0$. 
Hence we have 
\begin{align}\label{uaksurjk}
\left(([h^{i,i}]-[h^{j,j}])^2-\frac{1}{5}(6[h^{i,i}]-1)\right)*[E^{i,i}]=0
\end{align}
for any $1\leq i\leq d$. 

Let $W$ be an irreducible $A(\tripeven)$-module and $W_{\lambda_1,\ldots,\lambda_d}$ the simultaneous eigenspace for $[h^{i,i}]$ of eigenvalues $\lambda_i$ for all $1\leq i\leq d$. 
The identity \eqref{urwls} gives 
\begin{align}\label{tvcwpghre}
[h^{i,i}]^2*([h^{i,i}]-1)*(8[h^{i,i}]+1)*(8[h^{i,i}]-3)=0.
\end{align}
Hence we may assume that $\lambda_1,\ldots,\lambda_d\in\{0,-\frac{1}{8},1,\frac{3}{8}\}$. 

First we consider the case $\lambda_i=1$ or $\frac{3}{8}$ for some $i$. 
We then can assume that $\lambda_1=1$ or $\frac{3}{8}$ if necessary by permutating the pairs $(e^i,f^i)$ in the canonical basis $\{(e^i,f^i)\}_{1\leq i\leq d}$.
By the classification of irreducible $A(\Teven)$-modules, there exists $0\neq w\in W_{\lambda_1,\ldots,\lambda_d}$ such that $[E^{1,1}].w\neq 0$. 
It follows from \eqref{uaksurjk} that for any $j\geq 2$, the eigenvalue $\lambda_j$ must satisfy the equation 
\[
(\lambda_1-\lambda_j)^2=\frac{1}{5}(6\lambda_1-1).
\] 
This equation shows that $\lambda_j=0,-\frac{1}{8}$ if $\lambda_1=1,\frac{3}{8}$ respectively. 
Therefore, we may assume that $W_{\lambda,\mu,\ldots,\mu}\neq0$ for $\lambda= 1$ (resp. $\frac{3}{8}$) and $\mu=0$ (resp. $-\frac{1}{8}$). 

We now recall Remark \ref{remark4}.
By this remark, we see that there exist nonzero vectors $x^1,y^1\in W_{\lambda,\mu,\ldots,\mu}$ such that   
\[
v_{\lambda}^{1,1}(e^1,f^1).x^1=x^1,\quad v_{\lambda}^{2,1}(e^1,f^1).x^1=y^1,
\] 
where $v_{\lambda}^{k,l}(e^1,f^1) $ is the vector $v_{\lambda}^{k,l}$ defined in the proof of Proposition \ref{laisdqli} with $e=e^1,f=f^1$ for $k,l=1,2$.
Let $\alpha=\alpha_{\lambda}=\lambda-\mu$ and set   
\[
x^j:=\frac{1}{\alpha}[h^{j,1}].x^1,\quad y^j:=\frac{1}{\alpha}[h^{1,j}].y^1
\]
for $j\geq 2$.
\begin{lemma}\label{ienwsl2}
For $1\leq i,j,k\leq d$ with $i\neq j$, 
\begin{alignat*}{4}
[f^{i,j}].x^{k}&=\alpha(\delta_{j,k}y^{i}-\delta_{i,k}y^{j}),&\qquad [f^{i,j}].y^{k}&=0,\\
[e^{i,j}].x^{k}&=0,&\qquad [e^{i,j}].y^k&=-\alpha(\delta_{j,k}x^{i}-\delta_{i,k}x^{j}),\\
[h^{i,j}].x^{k}&=\delta_{j,k}\alpha x^{i},&\qquad [h^{i,j}].y^k&=-\delta_{i,k}\alpha y^{j}.
\end{alignat*} 
\end{lemma} 
\begin{proof}
We fix $i>1$. 
Since $y^1=v^{2,2}_\lambda(e^1,f^1).y^1$, Remark \ref{remark4} shows that $[F^{1,1}].y^1=0$.
Moreover, since $[h^{i,i}].y^1=\mu y^1$ for either $\mu=0$ or $-\frac{1}{8}$, we have $[F^{i,i}].y^1=0$.
On the other hand, \eqref{owklwqa} proves that
\begin{align*}
\alpha^2[f^{i,1}].y^{1}=([h^{1,1}]-[h^{i,i}])^2.[f^{i,1}].y^{1}=[h^{1,i}].[h^{i,1}].[f^{i,1}].y^{1}+[h^{i,1}].[h^{1,i}].[f^{i,1}].y^{1}. 
\end{align*}
Since $[h^{i,1}].[f^{i,1}].y^{1}=\frac{1}{2}[F^{1,1}].y^1=0$ and $[h^{1,i}].[f^{i,1}].y^{1}=-\frac{1}{2}[F^{i,i}].y^1=0$ by \eqref{ualru}, we have 
\begin{align}\label{mmes1}
[f^{i,1}].y^{1}=0.
\end{align}
If we exchange $y^1$ to $x^1$ in the argument above, then one has 
\[
[h^{i,1}].[f^{i,1}].x^{1}=\frac{1}{2}[F^{1,1}].x^1\quad\text{ and }[h^{1,i}].[f^{i,1}].x^{1}=[F^{i,i}].x^1=0.
\]
Thus by using \eqref{owklwqa} again, we get 
\begin{align*}
\alpha^2[f^{i,1}].x^{1}=[h^{1,i}].[h^{i,1}].[f^{i,1}].x^{1}+[h^{i,1}].[h^{1,i}].[f^{i,1}].x^{1}=\frac{1}{2}[h^{1,i}].[F^{1,1}].x^1. 
\end{align*}
Since $y^1=v_{\lambda}^{2,1}(e^1,f^1).x^1=\frac{1}{2\alpha^2}[F^{1,1}].x^1$ by Remark \ref{remark4}, we find that 
\begin{align}\label{mmes2}
[f^{i,1}].x^1=[h^{1,i}].y^1=\alpha y^i.
\end{align}
We can also prove that 
\begin{align}\label{mmes3}
[e^{i,1}].x^1&=0,\quad [e^{1,i}].y^{1}=\alpha x^i 
\end{align}
by the same method.
By using \eqref{pqkelo7} and the fact that $[f^{i,1}]^2=0$, we get 
\begin{align}
[f^{1,i}].x^i&=\frac{1}{\alpha}[f^{1,i}].[e^{1,i}].y^{1}=\alpha[y^1],\label{mmes4}\\
[f^{1,i}].y^i&=\frac{1}{\alpha}[f^{1,i}].[f^{i,1}].x^i=0.\label{mmes5}
\end{align}
As well we have 
\begin{align}\label{mmes6}
[e^{i,1}].x^i&=0,\quad [e^{i,1}].y^{i}=\alpha x^1.
\end{align}
Since $[h^{1,i}]*[e^{1,i}]=\frac{1}{2}[E^{1,1}]$ by \eqref{ualru}, we have 
\[
[h^{1,i}].x^i=\frac{1}{\alpha}[h^{1,i}].[e^{1,i}].y^1=\frac{1}{2\alpha}[E^{1,1}].y^1.
\]
Remak \ref{remark4} then proves that 
\begin{align}\label{mmes7}
[h^{1,i}].x^i=\alpha v^{1,2}_\lambda(e^1,f^1).y^1=\alpha x^1.
\end{align}
As well, we have 
\begin{align}\label{mmes8}
[h^{i,1}].y^i= \alpha y^1.
\end{align} 
These show that 
\begin{align}
[h^{1,i}].x^{1}=\frac{1}{\alpha}[h^{1,i}]^2.x^i=0,\label{mmes9}\\
[h^{i,1}].y^{1}=\frac{1}{\alpha}[h^{i,1}]^2.y^i=0.\label{mmes10}
\end{align} 

We calculate the identities 
\begin{align}
[h^{i,j}]*[h^{j,k}]&=[h^{i,k}]*([h^{k,k}]-[h^{j,j}]),\label{ttqwdd1}\\
[e^{i,j}]*[f^{j,k}]&= [h^{i,k}]*([h^{j,j}]-[h^{k,k}])\label{ttqwdd3}
\end{align}
which follow from \eqref{aliuhcq}, \eqref{uiahysdc1}--\eqref{uiahysdc6} for distinct integers $1\leq i,j,k\leq d$. 
Then we have 
\[
[f^{i,j}].x^j=\frac{1}{\alpha}[f^{i,j}].[e^{j,1}]y^1=[h^{1,i}].y^1=\alpha y^i
\]
for $i,j>1$ with $i\neq j$. 
As for $e^{i,j}$ and $h^{i,j}$, we get the desired identities in Lemma \ref{ienwsl2} by using the identities \eqref{mmes1}--\eqref{ttqwdd1} and \eqref{ttqwdd3}. 
Therefore, we see that Lemma \ref{ienwsl2} holds if $i=k$ or $j=k$. 

Finally let $1\leq i,j,k\leq d$ be mutually distinct. 
Then we see that 
\begin{align*}
\alpha[e^{i,j}].x^k=[e^{i,j}].[e^{j,k}].y^j=0,\quad \alpha[e^{i,j}].y^k=[e^{i,j}].[h^{j,k}].y^j=0,
\end{align*}
where we use that $[e^{i,j}].[e^{j,k}]=[e^{i,j}].[h^{j,k}]=0$. 
As well we have 
\begin{align*}
\alpha[f^{i,j}].x^k=[f^{i,k}].[h^{j,k}].x^k=0,\quad \alpha[f^{i,j}].y^k=[f^{i,j}].[f^{k,j}].x^k=0.
\end{align*}
By \eqref{ttqwdd3}, we see that 
\begin{align*}
\alpha[h^{i,j}].x^k&=[e^{i,k}].[f^{k,j}].x^k=-\alpha[e^{i,k}].y^j=0,\\
\alpha[h^{i,j}].y^k&=[e^{i,k}].[f^{k,j}].y^k=0.
\end{align*}
The proof is completed. 
\end{proof}
Compare Lemma \ref{ienwsl2} with \eqref{eedort1}--\eqref{eedort4}.  
In the case $\lambda=1$, we see that the linear map form $\Omega(\tripodd)$ to $\sum_{i=1}^d(\C x^i+\C y^i)$ defined by $e^i\mapsto x^i$ and $f^i\mapsto y^i$ is an $A(\tripeven)$-module homomorphism. 
Since the map is nonzero, it is an isomorphism.
In particular $W\cong\Omega(\tripodd)$. 
In the case $\lambda=\frac{3}{8}$, we see that $W\cong \Omega(\triptwodd)$.
Therefore, we have 
\begin{proposition}
Let $W$ be an irreducible $A(\tripeven)$-module.
If there exists a nonzero simultaneous eigenspace $W_{\lambda_1,\ldots,\lambda_d}$ for all $[h^{i,i}]$ such that $\lambda_i=1$ {\rm(}resp. $\frac{3}{8}${\rm)} for some $i$, then $W$ is isomorphic to $\Omega(\tripodd)$ {\rm(}resp. $\Omega(\triptwodd)${\rm)}.
\end{proposition}

We next consider the case that an irreducible $A(\tripeven)$-module $W$ includes a nonzero simultaneous eigenspace $W_{\lambda_1,\ldots,\lambda_d}$ such that $\lambda_i$ is either $0$ or $-\frac{1}{8}$ for any $1\leq i\leq d$. 
Then we have 
\begin{lemma}\label{pqwojqp}
For any $1\leq i,j\leq d$, $\lambda_i=\lambda_j$. 
Furthermore, $W=W_{\lambda,\ldots,\lambda}$ with $\lambda=\lambda_1$. 
\end{lemma}
\begin{proof}
Let $w\in W_{\lambda_1,\ldots,\lambda_d}$ be a nonzero vector. 
From the classification of irreducible $A(\Teven)$-modules, we see that 
\[
[E^{k,k}].w=[H^{k,k}].w=[F^{k,k}].w=0
\]
for any $1\leq k\leq d$. 
Let $1\leq i\neq j\leq d$. 
Then \eqref{ualru} proves  
\begin{align}\label{ooxrj2}
[h^{i,j}].[f^{i,j}].w=[h^{j,i}].[f^{i,j}].w=0.
\end{align}
Hence by \eqref{owklwqa}, 
\begin{align}\label{ttus}
0=[h^{i,j}].[h^{j,i}].[f^{i,j}].w+[h^{j,i}].[h^{i,j}].[f^{i,j}].w=(\lambda_i-\lambda_j)^2[f^{i,j}].w.
\end{align}
On the other hand \eqref{owq;autrf1} shows that 
\begin{align}\label{ooxjr}
2[e^{i,j}].[f^{i,j}].w=(\lambda_i-\lambda_j)^2w.
\end{align}
Therefore applying $[e^{i,j}]$ to the both sides in \eqref{ttus} proves $\lambda_i=\lambda_j$. 

Set $\lambda=\lambda_1$. 
By \eqref{pqkelo1}--\eqref{pqkelo6}, we see that $W_{\lambda,\ldots,\lambda}$ is closed under the actions of $[e^{i,j}]$, $[h^{i,j}]$ and $[f^{i,j}]$ for $1\leq i,j\leq d$. 
Thus by Proposition \ref{qloacql}, $W_{\lambda,\ldots,\lambda}$ is an $A(\tripeven)$-submodule of $W$, and $W=W_{\lambda,\ldots,\lambda}$.
\end{proof}

Now we consider the subspace $U:=[f^{i,j}].W$ for any $1\leq i\neq j<d$. 
We shall show that $U$ is an $A(\tripeven)$-submodule of $W$. 
Firstly, we have $[h^{k,k}]U\subset U$ for any $1\leq k\leq d$. 
Secondly, by \eqref{ooxrj2}, \eqref{ooxjr} and the fact that $[f^{i,j}]^2=0$, we see that $[x^{i,j}].U=0$ for any $x^{i,j}=e^{i,j},h^{i,j},f^{i,j}$.
It is clear that if $k,l\neq i,j$ then $[x^{k,l}].U\subset U$ for any $x^{k,l}=e^{k,l},h^{k,l},f^{k,l}$. 
Finally, for $1\leq k,l\leq d$ such that either $k$ or $l$ is in $\{i,j\}$, \eqref{ttqwdd3} implies that $[e^{k,l}].U=0$. 
Since $[f^{k,l}]*[f^{i,j}]=0$, we have $[f^{k,l}].U=0$. 
Furthermore we note that  
\[
[h^{k,i}].[f^{i,j}]=0,\quad [h^{i,k}].[f^{i,j}]=[f^{i,k}]*([h^{i,i}]-[h^{j,j}]).
\] 
Hence $[h^{k,l}].U\subset U$. 
Therefore we find that $[f^{i,j}].W$ is closed under the actions of $[e^{k,l}],[h^{k,l}],[f^{k,l}]$ for any $1\leq k,l\leq d$ and that it is an $A(\tripeven)$-submodule of $W$.  
Since $[f^{i,j}].W\neq W$ (otherwise $W=[f^{i,j}]^2W=0$), $[f^{i,j}].W$ is zero. 
Consequently we have $[f^{i,j}]=0$ on $W$. 
One can also prove that $[e^{i,j}]=[h^{i,j}]=0$ on $W$ for $1\leq i\neq j\leq d$. 
In particular, this and Proposition \ref{qloacql} show the following lemma. 
\begin{lemma}
The action of $A(\tripeven)$ on $W$ is commutative. 
Hence $W$ is one dimensional. 
\end{lemma}
Therefore, by taking a nonzero vector $w_\lambda\in W=W_{\lambda,\ldots,\lambda}$, we have an $A(\tripeven)$-module homomorphism from $\Omega(\tripeven)$ (resp. $\Omega(\triptweven)$) to $W$ defined by $\1\mapsto w_1$ (resp. $1_{\theta}\mapsto w_{-\frac{1}{8}}$). 
Consequently, we see that $W$ is isomorphic to either $\Omega(\tripeven)$ or $\Omega(\triptweven)$. 
\begin{theorem}\label{liawqwqq}
For the vertex operator algebra $\tripeven$ with $d\geq 2$, any irreducible $A(\tripeven)$-module is isomorphic to one in the list $\{\Omega(\trip^{\pm}),\Omega(\triptweo)\}$. 
\end{theorem}  

\section{Further structures of the vertex operator algebra $\tripeven$}
In this section we prove the irrationality of $\tripeven$ by constructing reducible  indecomposable $\tripeven$-modules, and determine the automorphism group of $\tripeven$. 
We also calculate the irreducible characters and their modular transformations. 
\subsection{Indecomposable $\tripeven$-modules}\label{Sectc4}
In this section we construct some reducible indecomposable $\Zpos$-gradable $\tripeven$-modules whose existence is shown in \cite{K1} in the case $d=1$. 
The existence of reducible indecomposable $\Zpos$-gradable $\tripeven$-modules proves that the vertex operator algebra $\tripeven$ is irrational.  

We recall the algebra $\mathcal{A}$ in Section \ref{Sectb1}. 
The vertex operator superalgebra $\trip$ is realized as the quotient $\mathcal{A}$-module $\mathcal{A}/\mathcal{A}_{\geq 0}$. 
We see that any $\trip$-module is an $\mathcal{A}$-module. 
Conversely, we can show that any $\mathcal{A}$-module $M$ is naturally an $\trip$-module if for any $u\in M$ and $\psi\in \h$, there is an integer $n_0$ such that $\psi_{(n)}u=0$ for $n\geq n_0$. 
In fact, for an $\mathcal{A}$-module $M$ satisfying this condition, the vertex operator \eqref{lqqljq} on $M$ is well defined.
We can then check that this gives an $\trip$-module structure on $M$. 

We consider the left ideal of $\mathcal{A}$ generated by $\psi(n)1$ for any vectors $\psi\in\h$ and $n\in\Zplus$, and denote it by $\mathcal{A}_{+}$. 
Then we see that for any $u\in\mathcal{A}$ and $\psi\in\h$, there exists $n_0\in\Zplus$ such that $\psi{(n)}u\in\mathcal{A}_{+}$ for any $n\geq n_0$. 
Therefore, the quotient $\mathcal{A}$-module  
\[
\indetrip:=\mathcal{A}/\mathcal{A}_+
\] 
becomes an $\trip$-module. 

The $\trip$-module $\indetrip$ is isomorphic to $\Lambda(\h\otimes\C[t^{-1}])$ as a vector space.
Set $\hat\1=1+\mathcal{A}_{+}$. 
Then $\indetrip$ is spanned by vectors of the form $\psi^1_{(-n_1)}\cdots\psi^r_{(-n_r)}\hat\1$ with $\psi^i\in\h$ and $n_i\in\Zpos$.
In particular, we have the decomposition 
\[
\indetrip=\bigoplus_{n=0}^{\infty}\indetrip_{(n)} 
\]
into a direct sum of generalized eigenspaces for $L_0$.
Hence $\indetrip$ is $\Zpos$-gradable. 
We note that the generalized eigenspace $\indetrip_{(0)}$ for $L_0$ of eigenvalue $0$ is spanned by $\psi^{1}_{(0)}\cdots\psi^r_{(0)}\hat{\1}$ for $\psi^i\in\h$ and $r\in\Zpos$.
Hence we may identify $\indetrip_{(0)}$ with the exterior algebra $\Lambda(\h)$.
We denote by ${}^r\Lambda(\h)$ the subspace spanned by vectors of the form $\psi^{1}_{(0)}\cdots\psi^s_{(0)}\hat{\1}$ for $\psi^i\in\h$ and $s\geq r$.
Consider the $\trip$-submodule $\indetrip[r]$ of $\indetrip$ generated from the subspace ${}^r\Lambda(\h)$ for any $r\in\Zpos$. 
It is clear that $\indetrip[r]\cong\trip\otimes{}^r\Lambda(\h)$ as vector spaces. 
Since $\dim\h=2d$, $\indetrip[2d+1]=0$. 
Thus we have a sequence of $\tripeven$-submodules
\[
0=\indetrip[2d+1]\subset \indetrip[2d]\subset\indetrip[2d-1]\subset\cdots\subset \indetrip[0]=\indetrip.
\]
By definition, for any $\psi\in\h$, $\psi_{(0)}\indetrip[r]\subset \indetrip[r+1]$.  
Thus $\psi_{(0)}$ acts trivially on the quotient $\indetrip[r]/\indetrip[r+1]$.
Therefore, 
\[
\indetrip[r]/\indetrip[r+1]\cong\trip\otimes ({}^r\Lambda(\h)/{}^{r+1}\Lambda(\h))
\]
as left $\trip$-modules. 
Since $\dim ({}^r\Lambda(\h)/{}^{r+1}\Lambda(\h))=\binom{2d}{r}$, $\indetrip[r]/\indetrip[r+1]$ is a direct sum of $\binom{2d}{r}$ copies of $\trip$ as an $\trip$-module.

We note that $\indetrip$ is naturally an $\tripeven$-module. 
We determine the space $\Omega(\indetrip)$ of singular vectors in $\indetrip$ as an $\tripeven$-module. 
If $u\in\Omega(\indetrip)\cap\indetrip[r]$ and $u\notin \indetrip[r+1]$ for some $r$, then $u+\indetrip[r+1]$ is a nonzero singular vector of $\indetrip[r]/\indetrip[r+1]$. 
Since this quotient is a direct sum of copies of $\tripeven$ and $\tripodd$, we see that $u+\indetrip[r+1]$ is a sum of eigenvectors for $L_0$ of weight $0$ or $1$.
This implies that $u\in\indetrip_{(0)}\oplus\indetrip_{(1)}$. 
Hence $\Omega(\indetrip)\subset \indetrip_{(0)}\oplus\indetrip_{(1)}$. 
 
It is clear that $\indetrip_{(0)}\subset \Omega(\indetrip)$. 
We claim that 
\begin{align}\label{ooe}
\Omega(\indetrip)\cap\indetrip_{(1)}=\indetrip[2d]_{(1)}.
\end{align} 
Let $v\in \indetrip_{(1)}$.
For a canonical basis $\{(e^i,f^i)\}$ of $\h$, we can find vectors $v_j\in \indetrip_{(0)}\,(1\leq j\leq 2d)$ such that $v=\sum_{i=1}^de^i_{(-1)}v_i+\sum_{j=1}^{d}f^{j}_{(-1)}v_{d+j}$. 
Then for any $\psi\in\h$ and $1\leq i\leq d$, 
\[
((e^i_{(-1)}\psi)_{(2)}+(e^i_{(-2)}\psi)_{(3)})v=-\psi_{(0)}v_{d+i},\quad ((f^i_{(-1)}\psi)_{(2)}+(f^i_{(-2)}\psi)_{(3)})v=\psi_{(0)}v_{i}.
\]
Therefore, if $v\in\Omega(\indetrip)$, then for each $i$, $v_i$ is annihilated by the action of $\psi_{(0)}$ for any $\psi\in\h$.
Thus $v_i\in{}^{2d}\Lambda(\h)\subset \indetrip[2d]$ and $\Omega(\indetrip)\cap\indetrip_{(1)}\subset \indetrip[2d]$. 
Since $\indetrip[2d]_{(1)}\subset\Omega(\indetrip)$, we get \eqref{ooe} and hence $\Omega(\indetrip)=\indetrip_{(0)}\oplus\indetrip[2d]_{(1)}$. 
\begin{proposition}\label{yiwkw}
The socle $\soc(\indetrip)$ of the $\tripeven$-module $\indetrip$ is $\indetrip[2d]$.  
\end{proposition}
\begin{proof}
Let $M$ be an irreducible $\tripeven$-submodule of $\indetrip$. 
Then $\Omega(M)\subset \Omega(\indetrip)$. 
Therefore, $M$ is generated from a vector in $\indetrip_{(0)}$ or $\indetrip[2d]_{(1)}$ by \eqref{ooe}. 
If $M$ is generated from a vector in $\indetrip[2d]_{(1)}$ then it is an irreducible submodule of $\indetrip[2d]$ isomorphic to $\tripodd$. 

Suppose that $M$ is generated from a vector $u$ in $\indetrip_{(0)}$. 
By the classification of irreducible $\tripeven$-modules, $M$ is isomorphic to $\tripeven$.
In particular, $M_{1}=0$. 
This implies that $(\psi_{(-2)}\phi)_{(1)}u=\phi_{(-1)}\psi_{(0)}u=0$ for any $\psi,\phi\in \h$. 
Therefore, we see that $u\in \indetrip[2d]_{(0)}$. 
Hence $M\subset\indetrip[2d]$. 
\end{proof}
The $\tripeven$-module $\indetrip$ is decomposable.
We see that the automorphism $\theta$ of $\mathcal{A}$ preserves the ideal $\mathcal{A}_+$. 
Hence $\theta$ acts on $\indetrip$.
If we denote by $\indetrip^\pm$ the $\pm1$-eigenspace of $\indetrip$ for $\theta$ respectively, then $\indetrip^\pm$ are $\tripeven$-modules and $\indetrip=\indetrip^+ \oplus\indetrip^-$ as $\tripeven$-modules. 
As a corollary of Proposition \ref{yiwkw}, we have 
\begin{corollary}
The $\tripeven$-modules $\indetrip^{\pm}$ are reducible and indecomposable.
\end{corollary}
\begin{proof}
The $\tripeven$-modules $\indetrip^\pm$ are reducible because $L_0$ does not act diagonally on them. 
We note that by Proposition \ref{yiwkw}, $\soc(\indetrip^{\pm})=\indetrip[2d]\cap\indetrip^{\pm}\cong\trip^{\pm}$ respectively. 
On the other hand, if $\indetrip^{+}$ or $\indetrip^{-}$ are decomposable then they must include at least two irreducible $\tripeven$-modules. 
This is a contradiction.    
\end{proof}

Therefore, the following is clear.
\begin{proposition}
The vertex operator algebra $\tripeven$ is irrational. 
\end{proposition}

\begin{remark}\label{wwes}
In the case $d=1$, $\Omega(\indetrip^+)$ gives an example of an $A(\Teven)$-module on which $[\w]$ does not act diagonally. 
Therefore, the ideal $A_0$ can not one dimensional. 
\end{remark} 

\subsection{Irreducible characters of $\tripeven$}\label{Secte2}
In this section we calculate characters of irreducible $\tripeven$-modules.
We also give their modular transformation.   

For a vertex operator algebra $V$, the character of an irreducible $V$-module $M$ with lowest weight $h$ is defined by 
\[
S_M(\tau)=\tr_{M} q^{L_0-\frac{c}{24}}=\sum_{n=0}^\infty \dim M_{h+n} q^{h-\frac{c}{24}+n},
\] 
where $q=e^{2\pi i\tau}$.
It is known that if the vertex operator algebra is $C_2$-cofinite then the character absolutely converges to a holomorphic function of $\tau$ on the upper half plane.
We denote the holomorphic function by $S_M(\tau)$ (see \cite{Zh1}).  
Let us consider the characters of irreducible $\tripeven$-modules $\trip^\pm$ and $\triptw^\pm$.  
It is easy to see that 
\begin{align*}
S_{\trip}(\tau)&=\left(q^{\frac{1}{24}}\prod_{n=1}^\infty(1+q^n)\right)^{2d}=\left(\frac{\eta(2\tau)}{\eta(\tau)}\right)^{2d},\\
S_{\triptw}(\tau)&=\left(q^{-\frac{1}{48}}\prod_{n=1}^\infty(1+q^{n-\frac{1}{2}})\right)^{2d}=\left(\frac{\eta(\tau)^2}{\eta(2\tau)\eta(\frac{\tau}{2})}\right)^{2d}
\end{align*}
with $q=e^{2\pi i\tau}$, where $\eta(\tau)=q^{\frac{1}{24}}\prod_{n=1}^\infty(1-q^n)$ is the Dedekind eta function. 
On the other hand, we have   
\begin{align*}
\tr_{\trip}\theta q^{L_0+\frac{d}{12}}&=\left(q^{\frac{1}{24}}\prod_{n=1}^\infty(1-q^n)\right)^{2d}=\eta(\tau)^{2d},\\
\tr_{\triptw}\theta q^{L_0+\frac{d}{12}}&=\left(q^{-\frac{1}{48}}\prod_{n=1}^\infty(1-q^{n-\frac{1}{2}})\right)^{2d}=\left(\frac{\eta(\frac{\tau}{2})}{\eta(\tau)}\right)^{2d}.
\end{align*}
We now set 
\begin{align}\label{modular}
\phi_1(\tau):=\frac{\eta(\tau)^2}{\eta(2\tau)\eta(\frac{\tau}{2})},\quad \phi_2(\tau):=\frac{\eta(\frac{\tau}{2})}{\eta(\tau)},\quad \phi_3(\tau)=\sqrt{2}\frac{\eta(2\tau)}{\eta(\tau)}
\end{align}
(see \cite[Chapter 4]{Wak}). 
Then we have 
\begin{align}
S_{\trip^\pm}(\tau)&=\frac{1}{2}\left(\frac{\phi_3(\tau)^{2d}}{2^d}\pm\eta(\tau)^{2d}\right),\label{iejkw}\\
S_{\triptw^\pm}(\tau)&=\frac{1}{2}\left(\phi_1(\tau)^{2d}\pm\phi_2(\tau)^{2d}\right).
\end{align} 
The modular transformations of the functions in \eqref{modular} are given by 
\begin{align}
\begin{split}\label{pqldjdf}
\phi_1(\tau+1)&=e^{-\frac{\pi i}{24}}\phi_2(\tau),\quad \phi_2(\tau+1)=e^{-\frac{\pi i}{24}}\phi_1(\tau),\quad \phi_3(\tau+1)=e^{\frac{\pi i}{12}}\phi_3(\tau),\\
\phi_1\left(-\frac{1}{\tau}\right)&=\phi_1(\tau),\quad \phi_2\left(-\frac{1}{\tau}\right)=\phi_3(\tau),\quad \phi_3\left(-\frac{1}{\tau}\right)=\phi_2(\tau)
\end{split}
\end{align}
which follow from the well known modular transformation lows 
\begin{align}\label{ir5ncfa}
\eta(\tau+1)=e^{\frac{\pi i}{12}}\eta(\tau),\quad \eta\left(-\frac{1}{\tau}\right)=(-i\tau)^{\frac{1}{2}}\eta(\tau).
\end{align}
By using the formula we have the following proposition.
\begin{proposition}\label{ilqlqao}
The modular transformations of $S_{\trip^\pm}(\tau)$ and $S_{\triptw^\pm}(\tau)$ with respect to the transformations $\tau\mapsto\tau+1$ and $\tau\mapsto-\frac{1}{\tau}$ are given by  
\begin{align*}
S_{\trip^\pm}(\tau+1)&=e^{\frac{d\pi i}{6}}S_{\trip^\pm}(\tau),\\
S_{\trip^\pm}\left(-\frac{1}{\tau}\right)&=\frac{1}{2^{d+1}}\left(S_{\triptweven}(\tau)-S_{\triptwodd}(\tau)\right) \pm\frac{(-i\tau)^d}{2}(S_{\tripeven}(\tau)-S_{\tripodd}(\tau)),\\
S_{\triptw^\pm}(\tau+1)&=\pm e^{-\frac{d\pi i}{12}}S_{\triptw^\pm}(\tau),\\
S_{\triptw^\pm}\left(-\frac{1}{\tau}\right)&=\frac{1}{2}\left(S_{\triptweven}(\tau)+S_{\triptwodd}(\tau)\right)\pm2^{d-1}\left(S_{\tripeven}(\tau)+S_{\tripodd}(\tau)\right). 
\end{align*}
\end{proposition}
\begin{proof}
We see that 
\[
\phi_1(\tau)^{2d}=S_{\triptweven}(\tau)+S_{\triptweven}(\tau),\quad\phi_2(\tau)^{2d}=S_{\triptweven}(\tau)-S_{\triptwodd}(\tau)
\]
and 
\[
\phi_3(\tau)^{2d}=2^{d}\left(S_{\tripeven}(\tau)+S_{\tripodd}(\tau)\right),\quad \eta(\tau)^{2d}=S_{\tripeven}(\tau)-S_{\tripodd}(\tau).
\]
Hence \eqref{iejkw}--\eqref{ir5ncfa} prove the proposition. 
\end{proof}
 
\subsection{The automorphism group of $\tripeven$}
We determine the automorphism group of $\tripeven$ in this section. 

We first recall that the group of all linear isomorphisms of $\h$ which preserve the skew-symmetric bilinear form $\langle\,\cdot\,,\cdot\,\rangle$ is the symplectic group $Sp(2d,\C)$.   
We extend the action of $Sp(2d,\C)$ on $\h$ to $\trip$ by the properties 
\begin{align*}
g(\1)&=\1,\\
g(\psi^1_{(-n_1)}\cdots\psi^r_{(-n_r)}\1)&=(g(\psi^1))_{(-n_1)}\cdots(g(\psi^r))_{(-n_r)}\1
\end{align*}
for any $g\in Sp(2d,\C)$, $\psi^i\in\h$ and $n_i\in\Zplus$. 
Then we have $gY(a,z)g^{-1}=Y(g(a),z)$ for any $g\in Sp(2d,\C)$ and $a\in\trip$. 
Since the definition of $\w$ does not depend on a choice of a canonical basis, we have $g(\w)=\w$ for any $g\in Sp(2d,\C)$.  
Therefore, any element of $Sp(2d,\C)$ induce an automorphism of $\trip$. 
In fact, the action is faithful, hence the automorphism group $\Aut(\trip)$ contains a subgroup isomorphic to $Sp(2d,\C)$. 
Conversely, we have 
\[
\langle g(\psi),g(\phi)\rangle\1=g(\psi)_{(1)}g(\phi)=g(\psi_{(1)}\phi)=\langle \psi,\phi\rangle\1
\] 
for any $g\in\Aut(\trip)$ and $\psi,\phi\in \h$. 
Hence any elements in $\Aut(\trip)$ give elements of $Sp(2d,\C)$.
This shows that 
\[
\Aut(\trip)\cong Sp(2d,\C).
\]   
We note that the automorphism $\theta$ is in the center of $Sp(2d,\C)$
and $\langle\theta\rangle$ is the center of $Sp(2d,\C)$.  
Therefore, $Sp(2d,\C)/\langle \theta\rangle$ faithfully acts on $\tripeven$.
We shall prove that $\Aut(\tripeven)\cong Sp(2d,\C)/\langle\theta\rangle$. 

We see that the characters $S_M(\tau)$ for $M=\trip^\pm,\triptw^\pm$ are mutually distinct. 
This implies that for any $g\in \Aut(\tripeven)$ and irreducible $\tripeven$-module $M$, the $\tripeven$-module $(M^g, Y^g(\,\cdot\,,z))$ with $M^g=M$ and $Y^g(\,\cdot\,,z)=Y(g(\cdot)\,,z)$ is isomorphic to itself because $S_{M^g}(\tau)=S_M(\tau)$. 
In particular, for any $g\in\Aut(\tripeven)$, there exists a unique $\tripeven$-module isomorphism $f_g:\tripodd\rightarrow(\tripodd)^g$ up to nonzero scalar multiple. 
Actually, if $f_g'$ is another $\tripeven$-module isomorphism from $\tripodd$ to $(\tripodd)^g$  then $f_g^{-1}\circ f_g'$ is in $\Hom_{\tripeven}(\tripodd,\tripodd)\cong\C$. 

We now consider the bilinear form $(\,\cdot\,,\cdot\,)_g$ on $\tripodd$ defined by $(u,v)_g=(f_g(u),f_g(v))$.
Then we see that it is a nondegenerate skew-symmetric, invariant bilinear form. 
Since $\tripodd$ is irreducible we have the following proposition (cf. \cite{Xu}). 
\begin{proposition}\label{iiwqkd}
There exists a nonzero constant $\alpha_g\in\C$ such that $(\,\cdot\,,\cdot\,)_g=\alpha_g(\,\cdot\,,\cdot\,)$.  
\end{proposition} 
\begin{proof}
The bilinear forms $(\,\cdot\,,\cdot\,)$ and $(\,\cdot\,,\cdot\,)_g$ satisfy that $(\tripodd_m, \tripodd_n)=(\tripodd_m,\tripodd_n)_g=0$ if $m\neq n$. 
This implies that the linear maps $\gamma$ and $\gamma_g$ from $\tripodd$ to the contragredient $\tripeven$-module $(\tripodd)'=\bigoplus_{n=1}^\infty(\tripodd_n)^{*}\subset D(\tripodd)$ defined by $\gamma(u)=(u,\cdot\,)$ and $\gamma_g(u)=(u,\cdot\,)_g$ respectively are $\tripeven$-module isomorphisms.
Hence $\gamma^{-1}\circ\gamma_g=\alpha_g\id_{\tripodd}$ for some $\alpha_g\in\C-\{0\}$. 
This proves the proposition.  
\end{proof}
By Proposition \ref{iiwqkd}, we can assume that $(\,\cdot\,,\cdot\,)_g=(\,\cdot\,,\cdot\,)$ if necessary by multiplying a suitable scalar to $f_g$. 
Hence we have $\langle f_g(\psi),f_g(\phi)\rangle=\langle\psi,\phi\rangle$ for any $\psi,\phi\in\h(\cong\tripodd_1)$. 
Therefore, the restriction of $f_g$ to $\tripodd_1$ gives an element of $Sp(2d,\C)$. 
We need to show the following lemma. 
\begin{lemma}
Let $\psi,\phi\in\h$ and $m,n\in\Zplus$. 
Then $g(\psi_{(-m)}\phi_{(-n)}\1)=f_g(\psi)_{(-m)}f_{g}(\phi)_{(-n)}\1$.
\end{lemma} 
\begin{proof}
First we assume that $n=1$. 
For any $u\in \tripeven$, we see that 
\[
(g(u),f_g(\psi)_{(-m)}f_{g}(\phi))=-(f_g(\psi)_{(m)}g(u),f_{g}(\phi))=(-1)^{m}(g(u)_{(m)}f_{g}(\psi),f_{g}(\phi)),
\]  
where the last identity follows from the skew symmetry formula 
\[
a_{(n)}b=(-1)^{kl}\sum_{i=0}^{\infty}\frac{(-1)^{n+1+i}}{i!}L_{-1}^{i}b_{(n+i)}a
\]
for $a\in\trip^{\bar{k}},b\in \trip^{\bar{l}}\,(k,l=0,1)$, $n\in\Z$ and the fact that $L_1f_{g}(\phi)=f_g(L_1\phi)=0$. 
Since $g(u)_{(m)}f_g(\psi)=f_g(u_{(m)}\psi)$ and the actions of $f_g$ and $g$ preserve the bilinear form $(\,\cdot\,,\cdot\,)$, we have 
\[
(g(u),f_g(\psi)_{(-m)}f_{g}(\phi))=(-1)^{m}(u_{(m)}\psi,\phi)=(u,\psi_{(-m)}\phi)=(g(u),g(\psi_{(-m)}\phi)).
\] 
Hence we have the lemma for $m\in\Zplus $ and $n=1$. 
As for arbitrary $n\in\Zplus$, by using Lemma \ref{twpdcnt} we can reduce to the case $n=1$. 
For example, 
\begin{align*}
g(\psi_{(-2)}\phi_{(-2)}\1)&=L_{-1}g(\psi_{(-2)}\phi)-g(\psi_{(-1)}\phi_{(-3)}\1)\\
&=L_{-1}f_g(\psi)_{(-2)}f_g(\phi)-f_g(\psi)_{(-1)}f_g(\phi)_{(-3)}\1=f_g(\psi)_{(-2)}f_g(\phi)_{(-2)}\1.
\end{align*}
The proof is completed. 
\end{proof}
This lemma states that if we extend $f_g$ to an automorphism $\tilde{f}_g$ of $\trip$ in the canonical way, then $\tilde{f}_g(a)=g(a)$ for any $a\in {\mathcal{L}}^2\tripeven$. 
In particular $g=\tilde{f}_g$ on $\tripeven_2\oplus\tripeven_3$. 
Hence Proposition \ref{llsp} proves that $g=\tilde{f}_g$ on $\tripeven$. 
In other words, the natural group homomorphism $Sp(2d,\C)\rightarrow \Aut(\tripeven)$ is surjective. 
Therefore, we have the following theorem.
\begin{theorem}
The automorphism group of $\tripeven$ is isomorphic to $Sp(2d,\C)/\langle\theta\rangle$. 
\end{theorem} 
   
\bibliographystyle{amsplain}

\end{document}